\documentclass[10 pt]{amsart}

\usepackage{amsmath}
\usepackage{amsxtra}
\usepackage{amscd}
\usepackage{amsthm}
\usepackage{amsfonts}
\usepackage{amssymb}
\usepackage{eucal}
\usepackage[all]{xy}
\usepackage{graphicx}
\usepackage{comment}
\usepackage{tikz}
\usepackage{tikz-cd}
\usepackage{thmtools}
\usepackage{thm-restate}
\usetikzlibrary{matrix, calc, arrows}
\usepackage{latexsym,amsmath,amscd,amssymb,epsfig,verbatim}
\usepackage[margin=1in,letterpaper,portrait]{geometry}

\usepackage[colorinlistoftodos, textsize=tiny]{todonotes}
\RequirePackage{color}
\definecolor{myred}{rgb}{0.75,0,0}
\definecolor{mygreen}{rgb}{0,0.5,0}
\definecolor{myblue}{rgb}{0,0,0.65}
	
\usepackage{hyperref}
  \hypersetup{colorlinks=true,
 	urlcolor=black,
	citecolor=blue}
  
\theoremstyle{plain}
  \newtheorem{thm}{Theorem}[section]
  \newtheorem{prop}[thm]{Proposition}
  \newtheorem{lem}[thm]{Lemma}
  \newtheorem{slem}{Sublemma}[lem]
  \newtheorem{cor}[thm]{Corollary}
  \newtheorem{conj}{Conjecture}
\theoremstyle{definition}
  \newtheorem{defn}[thm]{Definition}
  \newtheorem{note}[thm]{Notation}
  \newtheorem*{note*}{Notation}
  \newtheorem{rem}[thm]{Remark}
  \newtheorem{eg}[thm]{Example}
  
  \newtheorem{alg}[thm]{Algorithm}

\numberwithin{equation}{section}

\newcommand\nc{\newcommand}
\nc\on{\operatorname}
\nc\renc{\renewcommand}
\newcommand\ssec{\subsection}
\newcommand\sssec{\subsubsection}

\newcommand\BN{{\mathbb N}}

\newcommand\BF{{\mathbb F}}

\newcommand\BBZ{{\mathbb Z}}

\newcommand\BZ{{\mathbb Z}}

\newcommand{\id}{\mathrm{id}}

\newcommand\sminus{\text{-}}
\newcommand\dminus{\text{\--\,\--}}
\newcommand\splus{\text{+}}
\newcommand\dplus{\text{++}}
\newcommand\leftswap{\text{left}}
\newcommand\signswap{\text{sign}}
\newcommand\twoswap{\text{two}}
\DeclareMathOperator\block{Block}
\DeclareMathOperator\eblock{EvenBlock}
\DeclareMathOperator\oblock{OddBlock}

\newcommand\pblock{\mathrm{Block}_+}
\newcommand\mblock{\mathrm{Block}_-}
\newcommand\epblock{\mathrm{EvenBlock}_+}
\newcommand\opblock{\mathrm{OddBlock}_+}
\newcommand\emblock{\mathrm{EvenBlock}_-}
\newcommand\omblock{\mathrm{OddBlock}_-}
\DeclareMathOperator\plussetblock{\textit B}
\DeclareMathOperator\minussetblock{\textit B}
\DeclareMathOperator\involution{minSwap}
\def\Swap{\operatorname{Swap}}
\def\swap{\operatorname{swap}}
\def\sgn{\operatorname{sgn}}
\def\minus{\operatorname{minus}}

\def\listtodoname{List of Todos}
\def\listoftodos{\@starttoc{tdo}\listtodoname}

\title{Proof of Stasinski and Voll's Hyperoctahedral Group Conjecture}
\author{Aaron Landesman}

\AtEndDocument{\bigskip{\footnotesize%
  \textsc{Department of Mathematics, Stanford University, Stanford, CA 94305, U.S.A.}\par
  \textit{E-mail address}, \texttt{aaronlandesman@gmail.com}
}}

\usepackage{amsmath}
\begin{document}
\begin{abstract}
In a recent paper, Stasinski and Voll introduced a length-like statistic on hyperoctahedral groups and conjectured a product formula for this statistic's signed distribution over arbitrary quotients. Stasinski and Voll proved this conjecture for a few special types of quotients. We prove this conjecture in full, showing it holds for all quotients. In the case of signed permutations
with at most one descent, this formula gives the Poincar\'e polynomials for the varieties of symmetric matrices of a fixed rank.
\end{abstract}
\maketitle

\section{Introduction}
In this paper, we prove a conjecture of Stasinski and Voll \cite[Conjecture 1]{stasinski} as written in Theorem ~\ref{result}.
This result has relations to the Poincar\'e polynomials of the $n\times n$ symmetric matrices over $\BF_q$ of a given rank, as noted in \cite[p.\ 3]{stasinski}. The result also gives a formula for representation zeta functions associated to a particular group scheme, as described in \cite[p.\ 3-4]{stasinski} and \cite[Definition 1.2, Theorem C]{stasinski_2}. Additionally, this result immediately implies the validity of an interesting identity given in \cite[Proposition 5.5]{stasinski_2}. 

We now recall some notation, extrapolated upon in Section~\ref{background}, in order to state Theorem~\ref{result}.
Let $\mathbb N_0$ denote the non-negative integers, let $[n]_0$ denote $\left\{ 0,\ldots, n \right\}$, and let $\left[ \pm n \right]_0$ denote $\left\{ -n, \ldots, n \right\}$.
The {\it hyperoctahedral group} is by definition $B_n = \{w \in S_{2n+1} \mid  w(-i) = -w(i) \text{ for all } i \in [n]_0 \}.$
The key statistic on the hyperoctahedral group we will study is
the $L$ statistic, with $L: B_n \rightarrow \mathbb N_0$ defined by 
\begin{align*}
	L(w) = \frac{1}{2}\left| \left\{(i,j)\in [\pm n]_0 \times [\pm n]_0  \mid  i<j,w(i)>w(j),i \not\equiv j \bmod 2 \right\} \right|.
\end{align*}
Note that $L(w)$ is always an integer because $(i,j)$ satisfies $i<j,w(i)>w(j),i \not\equiv j \bmod 2$ if and only if $(-j,-i)$ does.

Intuitively, the statistic $L(w)$ measures how far $w$ is from the identity element, when consideration is restricted to indices of opposite parity.
For comparison, it is somewhat similar to the Coxeter {\em length function}
$l:B_n \rightarrow \BN_0$ which, for $w \in B_n$, which has the explicit formula 
\begin{align*}
	l(w) = \left | \left\{ \left( i,j \right) \in \left[ n \right]\times \left[ n \right] \mid i < j, w(i) > w(j) \right\} \right | 
	+ \left | \left\{ \left( i,j \right) \in \left[ n \right]\times \left[ n \right] | i \leq j, w(-i) > w(j) \right\} \right |,
\end{align*}
as shown in \cite[Proposition 8.1.1 and (8.2)]{bjorner}.
\begin{eg}
	As a quick example to understand the $L$ statistic, consider the element $w \in B_4$ given by $w = [2, -3, 1, 4]$ (so $w(1) = 2, w(2) = -3, w(3) = 1, w(4) = 4$). Then, $L(w) = 3$ corresponding to the $6$ pairs $(1,2),(-2,-1), (-2,1), (-1, 2), (-2,3),$ and $(-3,2)$.
\end{eg}

We next recall some standard notation associated to Coxeter groups. Let $D(w) =\{i \in [n-1]_0\mid l(ws_i)<l(w)\}$ denote the {\it descent set} of $w$,
where $s_i$ are the standard Coxeter generators for $B_n$ (see Lemma~\ref{lemma:coexeter-structure}).
For $I = \{i_1,\ldots,i_l\}_< \subset [n-1]_0,$
(where the $<$ subscript means $i_1 < \cdots< i_l$)
notate the {\it quotient} $B_n^I = \{w \in B_n \mid D(w) \subset I^c\},$ where $I^c$ denotes the complement of $I$ in $[n-1]_0$.

We also briefly recall the notation for $f_{n,I}(X)$ introduced in \cite[p.\ 2]{stasinski}. Define
$(\underline n)$ to be $1 - X^n$ if $n > 0$ and $1$ if $n = 0$. Further, define $(\underline n)! =(\underline n) \cdot (\underline {n-1}) \cdots (\underline 2)\cdot (\underline 1)$ and $(\underline n)!! = \prod_{i \in [n],i \equiv 0 \bmod 2} (\underline i).$
Given $I = \{i_1,\ldots,i_l\}_< \subset [n-1]_0,$ let $i_0 = 0,$ and $i_{l+1} = n.$ For $0 \leq t \leq l,$ define $\delta_t = i_{t+1} - i_t$
and notate
\begin{equation} 
\label{f_equation}
	f_{n,I}(X) = \frac{(\underline n)!}{(\underline {i_1})!\cdot \prod_{k=1}^{l}(\underline {\delta_k})!!}.
\end{equation}
With this notation in hand, we can now state our main result.
\begin{restatable}{thm}{mainresult}
\label{result}
Let $n \in \BN,I \subset [n-1]_0.$ Then
\begin{align}
	\label{equation:result}
	\sum_{w \in B_n^{I^c}}^{}(-1)^{l(w)}X^{L(w)} = f_{n,I}(X).
\end{align}
\end{restatable}

In \cite[Theorem 2]{stasinski}, Theorem ~\ref{result} was proved in the following special cases:
$I = \{0\}, I = [n-1]_0,$ and all elements of $I\cup \{n\}$ are even. The remainder of the paper presents a proof that holds for any $I \subset [n-1]_0$. The proof proceeds by induction on $n.$ We first demonstrate a recursive relation for the polynomials $f_{n,I}(X)$, and then show that this same recursion applies to the expression $\sum_{w \in B_n^{I^c}} (-1)^{l(w)}X^{L(w)}.$
A type-A analog of Theorem ~\ref{result}, where the sum is taken over elements of $S_n^{I^c}$ instead of $B_n^{I^c}$, can be found in \cite[Conjecture C]{klopsch}.

When this paper was in the final stages of preparation for submission, we learned that \cite[Conjecture C]{klopsch} and \cite[Conjecture 1]{stasinski} had recently been independently proven in a different way by Brenti and Carnevale \cite[Theorems 4.2 and 5.4]{brenti-carnevale}.
Their proof uses induction together with the additional technique of ``shifting and compressing,''
as described in \cite[Section 3]{brenti-carnevale}.

We now give an overview of our strategy for proving Theorem ~\ref{result}.
To start, we show the polynomials $f_{n,I}(X)$ satisfy a certain recurrence relation. 
The remainder of the paper is devoted to showing
the left hand side of ~\eqref{equation:result} satisfies the same recurrence relation.
In order to show this, we reduce the support of the left hand side to certain ``initially uncancelled'' and ``finally uncancelled'' elements.
To achieve this reduction, we introduce several involutions known as ``swaps'' which pair up permutations having opposite signs and the same value of
$L$, so that their contributions to the left hand side of ~\eqref{equation:result} cancel.
The initially uncancelled and finally uncancelled elements are characterized as those elements that have no swaps of certain types.
We then show there is a bijection between certain initially uncancelled elements of $B_n^{I^c}$ and certain finally uncancelled
elements of $B_{n-1}^{I^c}$, which is used to obtain the desired recurrence relation.

The outline of this paper is as follows. In Section ~\ref{background}, we present background material, including the relevant previous partial progress toward the main result, Theorem ~\ref{result}. In Section ~\ref{fni_recursion}, we give a recursion that the $f_{n,I}$ satisfy. Next, in Section ~\ref{combinatorial_induction} we state Theorem ~\ref{induction_result}, and show how this implies Theorem ~\ref{result}. In Section ~\ref{swapping}, we describe some crucial involutions which allow us to restrict the set of elements over which the generating function is summed. In Sections ~\ref{case_1},~\ref{case_2}, and ~\ref{case_3}, we show why the Parts 1, 2, and 3 Theorem ~\ref{induction_result} hold. By far the most difficult is Part 3, which involves characterizing certain elements in $B_n,$ which we call ``initially uncancelled'' and ``finally uncancelled,'' and then describing a bijection between these two types of elements. Finally, in Section ~\ref{further_conjecture}, we state a further conjecture, which involves a generalization of the generating function to two variables.

\section{Background}
\label{background}

\ssec{Notation for $B_n$}

Throughout this paper, we will only use congruence conditions of the form $a \equiv b \bmod 2,$ and so we drop the ``$\bmod \hspace{.07cm}2,$'' and henceforth write $a \equiv b$ to simplify notation. 
We let 
\begin{align*}
	[n] &=\{1,\ldots, n\}\\
	[n]_0 &= \{0,\ldots,n\}\\
	[\pm n]_0 &= \{-n,\ldots,n\}\\
	\BN &= \{n \in \BBZ,n > 0\} \\
	\BN_0 &= \BN \cup \{0\}.
\end{align*}

We have an action of $S_{2n+1}$ on $[\pm n]_0$ given by the standard permutation action on the $2n+1$ elements of $[\pm n]_0$.
Letting $B_n$ denote the hyperoctahedral group, as defined in the introduction, a convenient way to notate $w \in B_n$ will be $w = [i_1, \ldots, i_n]$ which means that for $t \in [n],$ we have $w(t) = i_t,w (0) = 0,$ and $w(-t) = -i_t.$

\begin{defn}
Let $w \in B_n.$ Define the {\it signed permutation matrix} which we notate as $w \in M_{n \times n}.$ It is defined by
\begin{align*}
	w_{i,j} = \begin{cases}
	1, &\text{ if }w(j) = i,\\
	-1 &\text{ if }w(j) = -i, \\
	0 &\text{ otherwise}.
\end{cases}
\end{align*}
\end{defn}

We often go back and forth between thinking about $w \in B_n$ and $w \in M_{n \times n}.$ When encountering the notation $w(k),$ think of $w \in B_n.$  
The following standard lemma describes the Coxeter structure of $B_n$.

\begin{lem}[\protect{\cite[Proposition 8.1.3]{bjorner}}]
	\label{lemma:coexeter-structure}
The group $B_n$ is the Coxeter group generated by the involutions $s_0,s_1,\ldots,s_{n-1}$ 
and satisfying the relations $s_i^2$ for $i \in [n-1]_0$, $(s_0s_1)^4,(s_is_{i+1})^3$ for $i \in [n-2],$ and $(s_is_j)^2$ if $j \neq i \pm 1.$
We may explicitly take $s_0 = [-1,2,3,\ldots,n],$ and $s_i = [1,2,\ldots, i-1,i+1,i,\ldots,n]$ for $i >0.$ 
\end{lem}

Fix $w \in B_n$
and let $l:B_n \rightarrow \BN_0$ denote the usual Coxeter length function. Define the sign function, $\sgn:B_n \rightarrow \{\pm 1\},$ by $\sgn(w) = (-1)^{l(w)}$.
It follows from this definition that $\sgn(w) = \det(w)$.
Throughout this paper, we will only have to deal with the easy-to-calculate $\sgn(w)$, instead of the slightly more difficult to calculate $l(w)$.
 Denote the {\it descent set} by $D(w) =\{i \in [n-1]_0\mid l(ws_i)<l(w)\},$ which is the usual Coxeter right descent set. Also, say $w$ {\it has a descent at $k$} if $k \in D(w).$ The descent set has an alternate characterization, which we now state in Lemma ~\ref{descent_equivalence}. We shall almost exclusively use the following characterization throughout this paper.

 \begin{lem}[\protect{\cite[Proposition 8.1.2]{bjorner}}]
\label{descent_equivalence}
For all $w \in B_n,$ there is an equality $D(w) = \{i \in [n-1]_0\mid w(i)>w(i+1)\}.$
\end{lem}

Let $k$ be a column index. Denote by $i_w(k)$ the unique row index such that $w_{i_w(k),k} \neq 0.$ That is, row $i_w(k)$ is the row in which the $\pm 1$ in column $k$ lies. Similarly, if $k$ is a row index, let $j_w(k)$ be the unique column index such that $w_{k,j_w(k)} \neq 0.$ When reference to some $w \in B_n$ is clear we often abuse notation by writing $i(k) = i_w(k)$ and $j(k) = j_w(k)$.

\begin{defn}
If $j$ and $l$ are column indices, say {\it $j$ is left of $l$} if $j < l$ and {\it $j$ is right of $l$} if $j > l$.
Further, for $s$ another column index, say {\it $s$ is column-between $j$ and $l$} when $\min(j,l) < s < \max(j,l)$.

Similarly, for row indices $i$ and $k$, say {\it $i$ is above $k$} if $i < k$
and {\it $i$ is below $k$} if $i > k$.
Finally, and for $t$ another row index, say {\it $t$ is row-between $i$ and $k$} when $\min(i,k) < t < \max(i,k)$
\end{defn}

\ssec{Relevant Results from Stasinski and Voll}

In order to proceed with the proof, we recount the relevant progress toward Theorem ~\ref{result} made in \cite{stasinski}. 
First, we recall that to prove Theorem ~\ref{result}, we can restrict to summing over chessboard elements, as defined in Definition ~\ref{chessboard_element}. Intuitively, if one colors the squares of a matrix $w \in B_n$ black and white, as on a chessboard, $w$ is a chessboard element if all nonzero elements occur on the white squares.

\begin{defn}
\label{chessboard_element}
For $n \in \BN$, define
$$C_n = \{w \in B_n \mid w_{i,j} \neq 0 \implies i+j \equiv 0\}.$$ 
If $w \in C_n,$ say $w$ is a {\it chessboard element}.
For $I \subset [n-1]_0,$ define $C_n^I = C_n \cap B_n^I$.
\end{defn}

In \cite{stasinski}, Stasinski and Voll call $C_n,$ defined in Definition ~\ref{chessboard_element}, $C_{n,0},$ in order to differentiate it from another set $C_{n,1} \subset B_n.$

\begin{lem}[Chessboard Lemma~\protect{\cite[p.\ 9, Lemma 8]{stasinski}}]
\label{chessboard}
For $n \in \BN$ and $I \subset [n-1]_0,$
\begin{align*}
	\sum_{w \in B_n^{I^c}}^{}(-1)^{l(w)}X^{L(w)} = \sum_{w \in C_n^{I^c}}^{}(-1)^{l(w)}X^{L(w)}.
\end{align*}
\end{lem}

\begin{defn}
For $n \in \BN$ and $I \subset [n-1]_0$ define $$S_{n,I}(X) = \sum_{w \in C_n^{I^c}}^{}(-1)^{l(w)}X^{L(w)}.$$
\end{defn}

As is immediate from the Chessboard Lemma ~\ref{chessboard}, we can restrict the sum in Theorem ~\ref{result} to chessboard elements:
\begin{cor}
\label{chessboard_sum}
To prove Theorem ~\ref{result}, it suffices to show $S_{n,I}(X) = f_{n,I}(X).$
\end{cor}

We next recall statistics introduced in ~\cite{stasinski}, which will help us calculate $L,$ using the upcoming Lemma ~\ref{abc_lemma}.

\begin{defn}
\label{abc_defn}	
For $n \in \BN, w \in B_n,$ and $j,j'$ so that $1 \leq j, j' \leq n,$ define the statistics
\begin{align*}
	a(w) &= |\{j \in [n] \mid w_{i(j),j}=-1,j\not\equiv 0\}|\\
	b_{j,j'}(w) &=
	\begin{cases}
	1, &\text{ if }j<j',i(j)>i(j'),j \not\equiv j'\\
	0, &\text{ otherwise}
\end{cases}\\
c_{j,j'}(w) &=
	\begin{cases}
	1, &\text{ if }j<j',w_{i(j'),j'}=-1,i(j)<i(j'),j \not\equiv j'\\
	0, &\text{ otherwise}
\end{cases}\\
b_{j'}(w)&= \sum_{j=1}^{n}b_{j,j'}(w) ,\\
c_{j'}(w) &= \sum_{j=1}^{n}c_{j,j'}(w),\\
 b(w)&=\sum_{j=1}^{n}b_j(w) ,\\
 c(w)&=\sum_{j=1}^{n}c_j(w).
\end{align*}
\end{defn}

\begin{rem}
It is immediate from the definitions that the statistics $a, b,$ and $c$ defined in Definition ~\ref{abc_defn} agree with those in \cite[Definition 4]{stasinski}.
Namely, 
\begin{align*}
	b(w) &= |\{(j,j') \in [n]\times [n] \mid j<j',i(j)>i(j'),j\not\equiv j'\}|\\
	c(w) &= |\{(j,j') \in [n]\times[n] \mid w_{i(j'),j'}=-1,j<j',i(j)<i(j'),j\not\equiv j'\}|.
\end{align*}
\end{rem}

\begin{lem}[abc Lemma~\protect{\cite[p.\ 7, Lemma 6, (11)]{stasinski}}]
\label{abc_lemma}
Let $w \in B_n.$ Then, $L(w) = a(w)+b(w)+2 c(w).$
\end{lem}

\section{The $f_{n,I}(X)$ recursion}
\label{fni_recursion}

Now that we have established preliminary notation and previous results, we will show the $f_{n,I}(X)$ satisfy a certain recursion. Following this section, we will show the polynomials $S_{n,I}(X)$ satisfy the same recursion and agree when $n = 1$, proving $S_{n,I}(X) = f_{n,I}(X).$ This then proves Theorem ~\ref{result} by Corollary ~\ref{chessboard_sum}.

\begin{defn}
\label{subtract_one_set}
Let $I = \{i_1,\ldots,i_l\}_< \subset [n-1]_0$. For $k \in [l+1]$, set
$$I^{(k)} = \{i_1,\ldots, i_{k-1},i_k-1,\ldots,i_l-1\} \cap [n-2]_0,$$
viewed as a subset of $[n-2]_0$.
\end{defn}

\begin{prop}
\label{recursion}
For $I = \{i_1,\ldots,i_l\}_< \subset [n-1]_0, t \in [l]_0$, recall $i_{l+1} = n, i_0 = 0$, and $\delta_t = i_{t+1} - i_t$, as defined prior to
~\eqref{f_equation}.
Let $m$ be the biggest index such that $\delta_m$ is odd, if such an $m$ exists. Otherwise, let $m = -1$. Then, the $f_{n,I}(X)$ satisfy the recurrence
\begin{align}
	\label{recursive_formula}
	f_{n,I}(X) &= - X^n \cdot f_{n-1,I^{(m+1)}}(X) +\sum_{t=m+1}^{l+1} X^{n-i_{t}} \cdot f_{n-1,I^{(t)}}(X).
\end{align}
\end{prop}
\begin{proof}
We shall only prove \eqref{recursive_formula} in the case that $m \neq -1$, as the proof of the case when $m = -1$ is analogous. Observe that by definition,
\begin{equation}
\label{double_fact_minus_one}
	(\underline {k-1})!! = \begin{cases}
		\frac{(\underline k)!!}{(\underline k)} &\text{ if }k\equiv 0,\\
		(\underline k)!! &\text{ if } k \not \equiv 0.
	\end{cases}
\end{equation}
For $t \in \BZ, 1 \leq t \leq l+1,$ 
	\begin{align*}
		f_{n-1,I^{(t)}}(X) &= \frac{(\underline {n-1})!}{(\underline {i_1})!\cdot (\underline {\delta_{t-1}-1})!!\prod_{k \in [l], k \neq t}(\underline {\delta_k})!!} & (\text{by } ~\eqref{f_equation})\\
	 &= \begin{cases}
		\frac{(\underline {n-1})! \cdot (\underline {\delta_{t-1}})}{(\underline {i_1})!\cdot \prod_{k=1}^{l}(\underline {\delta_k})!!} &\text{ if } \delta_{t-1} \equiv 0, \vspace{2mm} \\
		\frac{(\underline {n-1})!}{(\underline {i_1})!\cdot \prod_{k=1}^{l}(\underline {\delta_k})!!} &\text{ if  } \delta_{t-1} \not \equiv 0.
	\end{cases} & (\text{Using \eqref{double_fact_minus_one} with } k = \delta_{t-1}).
	\end{align*}
	In particular, since $\delta_m \not \equiv 0$ but $\delta_{t-1} \equiv 0,$ for $t-1 >m,$
	\begin{equation}
	\label{shifted_one}
		f_{n-1,I^{(t)}}(X) =
	\begin{cases}
		\frac{(\underline {n-1})! \cdot (\underline {\delta_{t-1}})}{(\underline {i_1})!\cdot \prod_{k=1}^{l}(\underline {\delta_k})!!} &\text{ if } t -1 > m, \vspace{2mm}\\
		\frac{(\underline {n-1})!}{(\underline {i_1})!\cdot \prod_{k=1}^{l}(\underline {\delta_k})!!} &\text{ if  } t -1 = m.
	\end{cases}
	\end{equation}
By substituting the right hand side of \eqref{shifted_one} into \eqref{recursive_formula}, we have to show
\begin{align}
	\frac{(\underline n)!}{(\underline {i_1})!\cdot \prod_{k=1}^{l}(\underline {\delta_k})!!} &= - X^n \cdot \frac{(\underline {n-1})! }{(\underline {i_1})!\cdot \prod_{k=1}^{l}(\underline {\delta_k})!!} + X^{n-i_{m+1}} \cdot \frac{(\underline {n-1})!}{(\underline {i_1})!\cdot \prod_{k=1}^{l}(\underline {\delta_k})!!}
\label{substituted_recursive_formula}	\\
	&\hspace{3 mm} + \sum_{t=m+2}^{l+1} X^{n-i_{t}} \cdot \frac{(\underline {n-1})! \cdot (\underline {\delta_{t-1}})}{(\underline {i_1})!\cdot \prod_{k=1}^{l}(\underline {\delta_k})!!}.
\nonumber
\end{align}
Dividing both sides of \eqref{substituted_recursive_formula} by $\frac{(\underline {n-1})! }{(\underline {i_1})!\cdot \prod_{k=1}^{l}(\underline {\delta_k})!!}$ implies that it is equivalent to show
\begin{align}
\label{divided_recursive_formula}
	(\underline n) &= X^{n-i_{m+1}}- X^n+\sum_{t=m+2}^{l+1} X^{n-i_{t}} \cdot (\underline {\delta_{t-1}}).
\end{align}
Recalling that $i_{l+1} = n$, we have a telescoping series:
\begin{align}
\label{telescoping}
	\sum_{t=m+2}^{l+1} X^{n-i_{t}} \cdot (\underline {\delta_{t-1}}) &= \sum_{t=m+2}^{l+1} X^{n-i_{t}} \cdot (1-X^{i_{t}-i_{t-1}})
	\\
	&=\sum_{t=m+2}^{l+1} \left(X^{n-i_{t}}-X^{n-i_{t-1}} \right)
	\nonumber\\
	&= \sum_{t=m+2}^{l+1} X^{n-i_t}-\sum_{t=m+1}^{l}X^{n-i_t}
	\nonumber\\
	&= X^{n-n} + \sum_{t=m+2}^{l} X^{n-i_t}-\sum_{t=m+2}^{l}X^{n-i_t} - X^{n-i_{m+1}}
	\nonumber\\
	& = 1 - X^{n-i_{m+1}}. \nonumber
\end{align}
Therefore, using \eqref{telescoping}, the right hand side of \eqref{divided_recursive_formula} is
\begin{align*}
	X^{n-i_{m+1}}- X^n+\sum_{t=m+2}^{l+1} X^{n-i_{t}} \cdot (\underline {\delta_{t-1}}) &= X^{n-i_{m+1}}- X^n  + 1 - X^{n-i_{m+1}} = 1 - X^n = (\underline n),
\end{align*}
implying \eqref{divided_recursive_formula} holds.
\end{proof}

\section{The Combinatorial Induction}
\label{combinatorial_induction}

To prove Theorem ~\ref{result}, that $f_{n,I}(X) = S_{n,I}(X)$, we will show that $S_{n,I}(X)$ satisfies the same recurrence established for $f_{n,I}(X)$ in Proposition ~\ref{recursion}. In this section, we will explain how Theorem ~\ref{induction_result} implies Theorem ~\ref{result}. This will be accomplished by partitioning $C_n^{I^c}$ according to where the nonzero entry in the bottom row is, and then analyzing the sums restricted to elements in this partition.
We first define notation needed to state Theorem ~\ref{induction_result}.

\begin{defn}
\label{cnik_defn}
Let $n \in \BN, I \subset [n-1]_0,$ with $I$ arbitrary, and $z \in [n]_0$. For convenience of notation, define $w_{n,0} = 0$. Define
\begin{align*}
	(C_n^{I^c})_z = \begin{cases}
	\{w \in C_n^{I^c} \mid w_{n,z}=1\}, &\text{ if }z\equiv n,\\
	\{w \in C_n^{I^c} \mid w_{n,z+1} = -1\}, &\text{ if } z \not \equiv n.
\end{cases}
\end{align*}
\end{defn}
Said another way, $(C_n^{I^c})_z$ are those elements of $C_n^{I^c}$ whose bottom nonzero entry lies either in column $z$ or $z+1$ (depending on the parity of $z$ and $n$). Therefore, 
$C_n^{I^c} = \cup_{z=0}^n (C_n^{I^c})_z$ defines a partition of $C_n^{I^c}$.
If $z < n$, by definition of $(C_n^{I^c})_z,$ for any $w \in (C_n^{I^c})_z$, we have $z \in D(w)$.

\begin{defn}
\label{snik_defn}
For $i_k \in \left\{ 0 \right\} \cup I \cup \left\{ n \right\}$, define 
\begin{align*}
S_{n,I,i_k}(X)=\sum_{w \in (C_n^{I^c})_{i_k}} (-1)^{l(w)}X^{L(w)}.
\end{align*}
\end{defn}

\begin{thm}
\label{induction_result}
Let $I = \{i_1,\ldots, i_l\}_<,i_0 = 0, i_{l+1} = n,$ and $\delta_k = i_{k+1} - i_k.$ Let $i_m$ be the greatest element of $I$ such that $\delta_m$ is odd, assuming such an $i_m \in I$ exists. If such an $i_m$ does not exist, let $m = -1$ and $i_m = -1.$ Then,
\begin{enumerate}
	\item $\sum_{i_k<i_m}^{} S_{n,I,i_k}(X)=0.$
	\item For $m \geq 0,$  $S_{n,I,i_m}(X)=-X^n \cdot f_{n-1,I^{(m+1)}}(X).$
	\item For $k \in \{m+1,\ldots, l+1\}$ we have $S_{n,I,i_k}(X)=X^{n-i_k} \cdot f_{n-1,I^{(k)}}(X).$
\end{enumerate}
\end{thm}

In the rest of this section, we show why Theorem ~\ref{induction_result} implies Theorem ~\ref{result}.

\begin{lem}
\label{bot_support}
We have the equality $S_{n,I}(X) = \sum_{i_k \in I \cup \{i_{l+1}\}}S_{n,I,i_k}(X).$

\end{lem}
\begin{proof}
	First, for each element $w,$ there must be exactly one $i \in [n]$ such that $w_{n,i} = \pm 1.$ Therefore, $S_{n,I}(X) = \sum_{z \in [n]_0}S_{n,I,z}(X).$ In fact, if $w \in C_n^{I^c},$ then $w \in \cup_{i_k \in I \cup \{i_{l+1}\}} (C_n^{I^c})_{i_k}$. This will prove the lemma, because $$S_{n,I}(X) = \sum_{z \in [n]_0}S_{n,I,z}(X)=\sum_{i_k \in I \cup \{i_{l+1}\}}S_{n,I,i_k}(X),$$
	as $S_{n,I,z} = 0$ if $z \notin I \cup \{i_{l+1}\}.$ 

Above, we are using that if $w \in C_n^{I^c}$ then $w \in (C_n^{I^c})_{i_k}$ for a unique $i_k \in I \cup \{i_{l+1}\}.$ This holds because there is
a unique nonzero entry in the bottom row of $w$.
\end{proof}

\begin{proof}
[Proof of Theorem ~\ref{result} assuming Theorem ~\ref{induction_result}] We proceed by induction on $n$. 

\vskip.1in
\noindent
{\sf The base case: Theorem ~\ref{result} holds when $n=1$.}

If $n = 1$, then either $I = \{0\}$ or $I = \emptyset$. If $n=1,I=\{0\},$ then Theorem ~\ref{result} holds by \cite[Theorem 2 Case 1]{stasinski} (though of course this is easy to alternatively compute by hand). If $n=1,I=\emptyset,$ Theorem ~\ref{result} holds because $C_n^{\emptyset^c} = \{\id \},$ so $S_{1,\emptyset}(X)=1=f_{1,\emptyset}(X)$.

\vskip.1in
\noindent
{\sf The inductive step: Theorem ~\ref{result} holds, assuming it holds with $n$ replaced by $n-1$.}

By Corollary ~\ref{chessboard_sum}, it suffices to prove that $S_{n,I}(X)$ satisfies the same recursion as $f_{n,I}(X).$ By Lemma ~\ref{bot_support}, 
$S_{n,I}(X) =\sum_{i_k \in I \cup \{i_{l+1}\}}S_{n,I,i_k}(X).$ We now have three cases, depending on whether $i_k < i_m, i_k = i_m,$ or $i_k > i_m$.

We will only prove the case that $m \geq 0$, as the case that $m = -1$ is analogous. Using Part 1 of Theorem ~\ref{induction_result}, it follows that for $i_t \in I$ with $i_t<i_m,$ we have $\sum_{i_t < i_m}^{} S_{n,I,i_t}(X) = 0.$ Combining this result with Lemma ~\ref{bot_support} gives us that $S_{n,I}(X) = \sum_{i_k \in I \cup \{i_{l+1}\},i_k \geq i_m}S_{n,I,i_k}(X).$ Next, using Part 2 of Theorem ~\ref{induction_result},
$$S_{n,I,i_m}(X) = -X^n \cdot f_{n-1,{I^{(m+1)}}}(X).$$
And finally, for $i_k > i_m,$ we use Part 3 of Theorem ~\ref{induction_result} to obtain  
$$S_{n,I,i_k}(X)=X^{n-i_k} \cdot f_{n-1,{I^{(k)}}} (X).$$
Therefore, by induction,
\begin{align*} 
S_{n,I}(X) &= -X^n \cdot f_{n-1,I^{(m+1)}}(X)+\sum_{k=m+1}^{l+1}X^{n-i_{k}} \cdot f_{n-1,I^{(k)}}(X)\\
&= -X^n \cdot S_{n-1,I^{(m+1)}}(X)+\sum_{k=m+1}^{l+1}X^{n-i_{k}} \cdot S_{n-1,I^{(k)}}(X).
\end{align*}
By Proposition ~\ref{recursion}, $S_{n,I}(X)$ and $f_{n,I}(X)$ satisfy the same recurrence. 
Hence, $f_{n,I}(X) = S_{n,I}(X)$, and Theorem ~\ref{result} holds by Corollary ~\ref{chessboard_sum}.
\end{proof}

\section{Swapping}
\label{swapping}
\ssec{Swapping Definitions and Examples}

We will prove Theorem \ref{induction_result} in Corollary ~\ref{part_1_induction_result}, Lemma ~\ref{part_2_induction_result}, and Corollary ~\ref{part_3_induction_result}. However, before doing so, we will need the Swapping Lemma, Lemma ~\ref{swapping_lemma}, to establish the properties of an involution on certain sets of ``swappable'' elements. 

Heuristically, swaps are pairs of columns, whose corresponding rows can be interchanged (swapped). This involution will preserve $L$ and change the sign of $w.$ Hence, we can use this involution to cancel the contribution of all ``swappable'' elements to the sum $S_{n,I}(X),$ and restrict the sum to a much smaller set.
When defining swaps, we further keep track of a column $k$ (and the resulting swaps are called $k$-swaps).
The column $k$ will often, though not always, be a column with $i(k) = n$ and $w_{n, k}= 1$, 
i.e., a column with a $1$ in its bottom row.
The main purpose of introducing $k$ is to relate elements in $B_n$ with a $1$ in column $k$ to elements of $B_{n-1}$, via
a somewhat involved algorithm discussed in Subsection ~\ref{subsection:lowering-map}, useful for our proof of Theorem ~\ref{result} 
by induction.

We will have several different types of swaps depending on the relative positions of the columns being swapped to column $k$ and
the signs of the entries in those columns. Sign swaps change the sign of a single entry in some column $\leq k$, left swaps interchange
two rows whose nonzero entries lie in columns $\leq k$, single-plus swaps are swaps involving a single column the right of $k$ that contains a $1$,
single-minus swaps are swaps involving a single column right of $k$ that contains a $-1$,
double-plus swaps are swaps involving two columns right of $k$ with $1$s, and double-minus swaps are swaps involving
two columns right of $k$ with $-1$s.
For the remainder of Section ~\ref{swapping}, fix $n \in \BN, I \subset [n-1]_0,w \in C_n,$ and $k \in [n]$ a column index.

\begin{defn}
\label{swap_defn}
Let $w \in C_n.$ A {\it $k$-swap for $w$} is a pair $(b,t) \in [n]_0^2$ such that all of the following hold:
\begin{enumerate}
	\item $b < t$,
	\item $b \equiv t$,
	\item if $s$ is a row index which is row-between $i(b)$ and $i(t)$, then $j(s) > k$,
	\item if $t > k,$ then $b \neq 0$,
	\item if $t > k$ and if $s$ is row-between $i(b)$ and $i(t)$, then $w_{s,j(s)} = -w_{i(t),t}$,
	\item if $b > k,$ then $w_{i(b),b} = w_{i(t),t}.$
\end{enumerate}
Let the set of {\it $k$-swaps for $w$} be $\Swap^k(w) = \{(b,t)\mid (b,t) \text{ is a $k$-swap for }w\}.$

We partition $\Swap^k(w)$ into the following named subsets, depending primarily on the relative positions of $b,t,$ and $k$ and the sign of $w_{i(t),t}$:
\begin{itemize}
	\item Let the set of {\it $k$-sign swaps for $w$} be $\Swap^k_{\signswap}(w) = \{(b,t) \in \Swap^k(w) \mid b = 0 \}.$
	\item Let the set of {\it $k$-left swaps for $w$} be $\Swap^k_{\leftswap}(w) = \{(b,t) \in \Swap^k(w) \mid b \neq 0, t \leq k \}.$
	\item Let the set of {\it $k$-single-plus swaps for $w$} be $\Swap^k_{\splus}(w) = \{(b,t) \in \Swap^k(w) \mid b \leq k < t, w_{i(t),t} = 1\}.$
	\item Let the set of {\it $k$-single-minus swaps for $w$} be $\Swap^k_{\sminus}(w) = \{(b,t) \in \Swap^k(w) \mid b \leq k < t, w_{i(t),t} = -1\}.$
	\item Let the set of {\it $k$-double-plus swaps for $w$} be $\Swap^k_{\dplus}(w) = \{(b,t) \in \Swap^k(w) \mid k < b, w_{i(t),t} = 1 \}.$
	\item Let the set of {\it $k$-double-minus swaps for $w$} be $\Swap^k_{\dminus}(w) = \{(b,t) \in \Swap^k(w) \mid k < b, w_{i(t),t} = -1\}.$
\end{itemize}
For $(b,t) \in \Swap^k(w),$ let $\{x,y\}_< = \{i(b),i(t)\}$ and define 
$$\swap_{(b,t)}(w) = s_{y-1} s_{y-2} \cdots s_{x+1} s_x s_{x+1} \cdots s_{y-1}w.$$
\end{defn}
\begin{rem}
\label{swap_explanation}
Note that when $b = 0$, $\swap_{(b,t)}(w)$ is almost the same matrix as $w,$ but with the sign of row $i(t)$ reversed. When $b \neq 0,\swap_{(b,t)}(w)$ is almost the same matrix as $w$, but with rows $i(b)$ and $i(t)$ interchanged. See Example ~\ref{swap_eg} for pictorial examples.
\end{rem}
\begin{eg}
\label{swap_eg}
We now give examples of various elements in the set of $k$-swaps. 
\begin{itemize}
\item
Let $w$ be as pictured below. Then, $(0,2) \in \Swap_{\signswap}^3(w)$ because $2\equiv 0$ and there are no entries in columns $\leq 3$ which lie above row $i(2) = 2.$
\[w=
	\begin{tikzpicture}[baseline=-0.5ex]
\matrix(A)[matrix of math nodes,column sep={20pt,between origins},row
    sep={20pt,between origins},left delimiter={(},right delimiter={)}] at (0,0)
  {
    & |[name = tswap]|& |[name = k]|& & 1\\
    &|[name=t]| -1 & & & \\
    1 & & & &\\
    & & &-1 & \\
     & & 1& & \\
  };
  \draw[->]
  ;
  \draw[->,gray];
  \node [yshift = .7 cm]at (k.north){$k$};
   \node [yshift = .7 cm]at (tswap.north){$t$};
   \end{tikzpicture}
\hspace{1cm} \swap_{(0,2)}(w)=
	\begin{tikzpicture}[baseline=-0.5ex]
\matrix(B)[matrix of math nodes,column sep={20pt,between origins},row
    sep={20pt,between origins},left delimiter={(},right delimiter={)}] at (6,0)
  {
    & |[name = tswap]|&|[name = k]| & & 1\\
    &|[name=t]| 1 & & & \\
    1 & & & &\\
    & & &-1 & \\
     & & 1& & \\
  };
  \draw[->] (t) circle (2mm)
  ;
  \draw[->,gray]
  ;\node [yshift = .7 cm]at (k.north){$k$};
  \node [yshift = .7 cm]at (tswap.north){$t$};
  \end{tikzpicture}
  .\]
\item
Let $w$ be as pictured below. Then, $(1,3)\in \Swap^4_{\leftswap}(w)$ because $1 \equiv 3,$ the only row-between row $i(1) = 3$ and row $i(3) = 5$ is row $4,$ and $j(4)=6$ is right of $4.$
\[w=
	\begin{tikzpicture}[baseline=-0.5ex]
\matrix[matrix of math nodes,column sep={20pt,between origins},row
    sep={20pt,between origins},left delimiter={(},right delimiter={)}] at (0,0)
  {
    |[name = bswap]|& &|[name = tswap]| & |[name = k]|& -1 &\\
    &1 & & & & \\
    |[name=bhigh]|1& & |[name=thigh]|& & &\\
     & & & & &1 \\
    |[name=blow]|& & |[name=tlow]|-1 & & & \\
    & & & 1& & \\
  };
  \node [yshift = .7 cm]at (k.north){$k$};
     \node [yshift = .7 cm]at (tswap.north){$t$};
        \node [yshift = .7 cm]at (bswap.north){$b$};
   \end{tikzpicture}
\hspace{1cm} \swap_{(1,3)}(w)=
	\begin{tikzpicture}[baseline=-0.5ex]
\matrix(B)[matrix of math nodes,column sep={20pt,between origins},row
    sep={20pt,between origins},left delimiter={(},right delimiter={)}] at (6,0)
   {
    |[name = bswap]|& & |[name = tswap]|&|[name = k]| & -1 &\\
    &1 & & & & \\
    |[name=bhigh]|\bullet& & |[name=thigh]|-1& & &\\
     & & & & &1 \\
    |[name=blow]|1& & |[name=tlow]| \bullet& & & \\
    & & & 1& & \\
  };
  \draw[->]
  (bhigh) to (blow);
  \draw[->]
  (tlow) to (thigh)
  ;
  \draw[->,gray]
  ;\node [yshift = .7 cm]at (k.north){$k$};
       \node [yshift = .7 cm]at (tswap.north){$t$};
        \node [yshift = .7 cm]at (bswap.north){$b$};
  \end{tikzpicture}.
  \]
\item Let $w$ be as pictured below. Then, $(1,3)\in \Swap^2_{\sminus}(w)$ because $1 \equiv 3$ and the row-between row $i(1) = 3$ and row $i(3) = 1$ is row $2$ with $w_{2,j(2)} = 1.$
\[w=
	\begin{tikzpicture}[baseline=-0.5ex]
\matrix[matrix of math nodes,column sep={20pt,between origins},row
    sep={20pt,between origins},left delimiter={(},right delimiter={)}] at (0,0)
  {
   |[name=bhigh]|&|[name = k]| & |[name=thigh]|-1& & &\\
    & & & 1& & \\
    |[name=blow]|1& &|[name=tlow]| & & &\\
     & & & & &1 \\
    & & & & 1& \\
    & 1& & & & \\
  };
  \draw[->,gray]
  ;\node [yshift = .7 cm]at (k.north){$k$};
       \node [yshift = .45 cm]at (thigh.north){$t$};
        \node [yshift = .7 cm]at (bhigh.north){$b$};
   \end{tikzpicture}
\hspace{1cm} \swap_{(2,6)}(w)=
	\begin{tikzpicture}[baseline=-0.5ex]
\matrix(B)[matrix of math nodes,column sep={20pt,between origins},row
    sep={20pt,between origins},left delimiter={(},right delimiter={)}] at (6,0)
    {
   |[name=bhigh]|1 &|[name = k]| & |[name=thigh]|\bullet& & &\\
    & & & 1& & \\
    |[name=blow]|\bullet& &|[name=tlow]| -1& & &\\
     & & & & &1 \\
    & & & & 1& \\
    & 1& & & & \\
  };
  \draw[->]
  (blow)to(bhigh);
  \draw[->]
  (thigh) to (tlow)
  ;
  \draw[->,gray]
  ;\node [yshift = .7 cm]at (k.north){$k$};
       \node [yshift = .52 cm]at (thigh.north){$t$};
        \node [yshift = .5 cm]at (bhigh.north){$b$};
  \end{tikzpicture}.
  \]
\item Let $w$ be as pictured below. Then, $(2,6)\in \Swap^1_{\dminus}(w).$ Note that $2 \equiv 6$ and the only rows between row $i(2) = 2$ and row $i(6) = 6$ are rows $3,4,$ and $5$, and $j(3), j(4)$ and $j(5)$ all lie right of $1$ and have $1$s in them.
\[w=
	\begin{tikzpicture}[baseline=-0.5ex]
\matrix[matrix of math nodes,column sep={20pt,between origins},row
    sep={20pt,between origins},left delimiter={(},right delimiter={)}] at (0,0)
  {
    |[name = k]|1&|[name = bswap]| & & & &|[name = tswap]|\\
    &|[name=bhigh]|-1 & & & & |[name=thigh]|\\
    & & & &1 &\\
     & & & 1& & \\
    & & 1 & & & \\
    & |[name=blow]|& & & & |[name=tlow]|-1\\
  };
  \node [yshift = .5 cm]at (k.north){$k$};
       \node [yshift = .7 cm]at (tswap.north){$t$};
        \node [yshift = .7 cm]at (bswap.north){$b$};
   \end{tikzpicture}
\hspace{1cm} \swap_{(2,6)}(w)=
	\begin{tikzpicture}[baseline=-0.5ex]
\matrix(B)[matrix of math nodes,column sep={20pt,between origins},row
    sep={20pt,between origins},left delimiter={(},right delimiter={)}] at (5.8,0)
   {
    |[name = k]|1&|[name = bswap]|& & & &|[name = tswap]|\\
    &|[name=bhigh]| \bullet& & & & |[name=thigh]|-1\\
    & & & &1 &\\
     & & & 1& & \\
    & & 1 & & & \\
    & |[name=blow]|-1& & & & |[name=tlow]|\bullet\\
  };
  \draw[->]
  (bhigh) to (blow);
  \draw[->]
  (tlow) to (thigh)
  ;
  \draw[->,gray]
  ;\node [yshift = .5 cm]at (k.north){$k$};
       \node [yshift = .7 cm]at (tswap.north){$t$};
        \node [yshift = .7 cm]at (bswap.north){$b$};
  \end{tikzpicture}
  \].
\end{itemize}
\end{eg}

\begin{rem}
The elements of $\Swap_{\dplus}^k(w)$ and $\Swap_{\dminus}^k(w)$ are reminiscent of type 1 ``odd sandwiches in $w$,'' as described in \cite[Definition 15]{stasinski}. Also, elements of $\Swap_{\signswap}^k(w)$ are reminiscent of type 2 ``odd sandwiches in $w$.'' We will soon define the maps \sloppy{$\involution^k_l$, $\involution^k_m$, and $\involution^k_p$} in Definition \ref{swapping_involution_defn}, which are analogs of the odd sandwich involution $w \mapsto w^\vee$ as defined in \cite[Definition 17]{stasinski}. In Subsections ~\ref{the_swapping_lemma} and ~\ref{the_swapping_involution}, we develop results for swaps similar to \cite[Lemma 19]{stasinski} for odd sandwiches.
\end{rem}

\begin{eg}
We now give an example of $\Swap^4(w)$ for $w = [-9,2,5,10,3,4,7,-6,1,-8].$ Writing $w$ in permutation matrix notation, 
$$w=\begin{pmatrix}
	 &  &  & &  &  &  &  & 1 & \\
	 & 1 &  &  &  &  &  &  &  & \\
	 &  &  &  & 1 &  &  &  &  & \\
	 &  &  &  &  & 1 &  &  &  & \\
	 &  & 1 &  &  &  &  &  &  & \\
	 &  &  &  &  &  &  & -1 &  & \\
	 &  &  &  &  &  & 1 &  &  & \\
	 &  &  &  &  &  &  &  &  & -1\\
	 -1&  &  &  &  &  &  &  &  & \\
	 &  &  & 1 &  &  &  &  &  & \\
\end{pmatrix}$$
Observe that $w \in C_{10},$ because $\forall i \in [10],i \equiv w(i).$ Also, $D(w)=\{0,4,7,9\}.$ 
We see 
\begin{itemize}
	\item $\Swap^4_{\signswap}(w)=\{(0,2)\}$,
	\item $\Swap^4_{\leftswap}(w)=\{(1,3)\}$,
	\item $\Swap^4_{\dminus}(w)=\{(8,10)\}$,
	\item $\Swap^4_{\splus}(w)=\{(1,7),(3,7)\}$,
	\item $\Swap^4_{\sminus}(w)=\Swap^4_{\dplus}(w) = \emptyset.$
\end{itemize}
As some non-examples, the following are not in the set of $4$-swaps for $w$:
\begin{itemize}
	\item $(2,3) \not\in \Swap^4_{\leftswap}(w),$ because $2 \not\equiv 3$,
	\item $(3,5) \not\in \Swap^4_{\splus}(w),$ because $w_{4,6}=1 \neq -w_{i(5),5}$,
	\item $(5,9) \not\in \Swap^4_{\dplus}(w)$ because row $2$ is row-between $i(5)=3$ and $i(9)=1$, but $j(2) = 2 \not > 4$.
\end{itemize}
\end{eg}

We now define a few more classes of swaps.
Informally, two-swaps are those involving two columns,
general-left swaps are those which only alter columns left of $k$ in $w$,
general-minus swaps are all swaps except single-plus swaps, and general-plus swaps are all swaps except single-minus swaps.

\begin{defn}
\label{general_swap_defn}
Define the set of {\it $k$-two-swaps for $w$} by 
$$\Swap^k_{\twoswap}(w) = \Swap^k(w) \setminus \Swap^k_{\signswap}(w).$$ 
Define the set of {\it general-left $k$-swaps for $w$} by
$$\Swap^k_l(w) = \Swap^k_{\signswap}(w) \cup \Swap^k_{\leftswap}(w).$$
Define the set of {\it general-minus $k$-swaps for $w$} by
$$\Swap^k_m(w)=\Swap^k(w) \setminus \Swap^k_{\splus}(w).$$
Define the set of {\it general-plus $k$-swaps for $w$} by
$$\Swap^k_p(w) = \Swap^k(w) \setminus \Swap^k_{\sminus}(w).$$
\end{defn}

\ssec{The Swapping Lemma}
\label{the_swapping_lemma}
Having defined various sets of swaps, we next show that if $(b,t) \in \Swap^k(w)$, then the two chessboard elements $w$ and $\swap_{(b,t)}(w)$ have descent sets differing by at most one element, have the same value of $L,$ and have opposite signs. We can then use these properties to cancel out certain elements $w$ with $\Swap^k(w) \neq \emptyset$ from $S_{n,I}(X).$

\begin{lem}
\label{descent_set_swap}
Let $w \in C_n.$ For any $(b,t) \in \Swap^k(w),$
we have $D(w) \setminus \{k\} = D(\swap_{(b,t)}(w)) \setminus \{k\}.$
\end{lem}
\begin{proof}
First, note that $i \in D(w)$ if and only if $w(i) > w(i+1).$ Since $w(i) = (\swap_{(b,t)}(w))(i)$ for $i \notin \{b,t\},$ we see that for $i \in [n-1]_0 \setminus \{b-1,b,t-1,t\},$ we have $i \in D(w) \iff i \in D(\swap_{(b,t)}(w)).$ It only remains to deal with $i \in \{b-1,b,t-1,t\}.$ Now, let $h \in \{b,t\},m \in \{h+1, h-1\},$ and $i = \min(h,m).$ The only way we can fail to have $i \in D(w) \iff i \in D(\swap_{(b,t)}(w))$ is when $w(m)$ takes a value between $w(h)$ and $(\swap_{(b,t)}(w))(h).$ If $h \leq k,$ then Property (3) in Definition ~\ref{swap_defn} ensures $w(m)$ is not between $w(h)$ and $(\swap_{(b,t)}(w))(h),$ unless $i = h = k.$ If $h > k,$ then Properties (5) and (6) of Definition ~\ref{swap_defn} ensure $w(m)$ is not between $w(h)$ and $(\swap_{(b,t)}(w))(h),$ unless $i = m = k.$ Therefore, if $i \neq k$, we have $i \in D(w) \iff i \in D(\swap_{(b,t)}(w)),$ so $D(w) \setminus \{k\} = D(\swap_{(b,t)}(w)) \setminus \{k\}.$
\end{proof}

\begin{lem}
\label{swap_after_swap}
Let $w \in C_n$. If $(b,t)\in \Swap^k(w)$, then $(b,t) \in \Swap^k(\swap_{(b,t)}(w)).$
\end{lem}
\begin{proof}

First, observe that if $(b,t) \in \Swap^k_\signswap(w),$ we have $b = 0$ and by Property (4) of Definition ~\ref{swap_defn}, $t \leq k$, so the last two conditions are vacuously satisfied. Since $\swap_{(b,t)}(w)$ only differs from $w$ in column $t,$ by Remark ~\ref{swap_explanation} the first three conditions of Definition ~\ref{swap_defn} hold, as they are independent of the value of $w_{i(t),t}.$ If $(b,t) \in \Swap^k_\twoswap(w)=\Swap^k(w) \setminus \Swap^k_\signswap(w)$, by Remark ~\ref{swap_explanation}, $i(b)$ and $i(t)$ are interchanged, and everything else is kept constant. Since all six conditions of Definition ~\ref{swap_defn} still hold after $i(b)$ and $i(t)$ are interchanged, $(b,t) \in \Swap^k(\swap_{(b,t)}(w)).$
\end{proof}

\begin{lem}[Swapping Lemma]
\label{swapping_lemma}
If $w \in C_n$ and $(\beta,\tau)\in \Swap^k(w),$ then\footnote{For this proof only, use $(\beta,\tau)$ instead of $(b,t),$ because we wish to reserve $b$ for the function $b$ as defined in Definition ~\ref{abc_defn} and we don't want to use the confusing subscript $b_b.$}
\begin{align*}
	L(w) &= L(\swap_{(\beta,\tau)}(w))\\
	\sgn(w) &= -\sgn(\swap_{(\beta,\tau)}(w)).
\end{align*}

\end{lem}
\sssec*{Idea of the proof of Lemma ~\ref{swapping_lemma}}
The most difficult part of this proof is showing $L(w) = L(\swap_{(\beta,\tau)}(w))$ when $(\beta,\tau)\in \Swap^k_\twoswap(w).$ To do this, the idea is to break up $L(w)$ and $L(\swap_{(\beta,\tau)}(w))$ as a sum of terms of the form $a(w), b_{j,j'}(w),$ and $c_{j,j'}(w).$ Since only columns $\beta,\tau$ differ between $w$ and $\swap_{(\beta,\tau)}(w)$ we will have $b_{j,j'}(w) = b_{j,j'}(\swap_{(\beta,\tau)}(w))$ and $c_{j,j'}(w) = c_{j,j'}(\swap_{(\beta,\tau)}(w))$ when $\{\beta,\tau\} \cap \{j,j'\} = \emptyset.$ Finally, by examining the relative position of $\beta,\tau, i(\beta)$ and $i(\tau)$, we will compute that the remaining terms in $L(w)$ and $L(\swap_{(\beta,\tau)}(w))$ sum to the same number.

\begin{proof}
We have $\sgn(w) = -\sgn(\swap_{(\beta,\tau)}(w)),$ since by Definition ~\ref{swap_defn}, $\swap_{(\beta,\tau)}(w)$ is the product of an odd number of transpositions with $w$.

So, we only have to check that $L(w) = L(\swap_{(\beta,\tau)}(w))$. For this, we use the abc Lemma ~\ref{abc_lemma} and compute how $a,b,$ and $c$ change. There are now two cases to verify, depending on whether $\beta = 0$. 

\vskip.1in
\noindent
{\sf The case of $(\beta,\tau) = (0,\tau) \in \Swap^k_\signswap(w)$.}

In this case, $a(w) = a(\swap_{(0,\tau)}(w))$, because $\tau \equiv 0$ by Property (2) of Definition ~\ref{swap_defn}, and $a$ only counts the number of $-1$s in odd columns. Next, $b(w) = b(\swap_{(0,\tau)}(w))$, because only the sign of row $i(\tau)$ changes. And finally, we check $c(w) = c(\swap_{(0,\tau)}(w)).$ By Property (4) of Definition ~\ref{swap_defn}, $k > \tau$. Hence, by Property (3) of Definition ~\ref{swap_defn}, there are no column indices $j \in [n]$, $j<\tau,i(j)<i(\tau),$ so $c_\tau(w) = c_\tau(\swap_{(0,\tau)}(w)) = 0.$ Additionally, it is evident that $c_{j'}(w) = c_{j'}(\swap_{(0,\tau)}(w))$ for $j' \neq \tau.$ Hence, $c(w) = c(\swap_{(0,\tau)}(w))$.
Therefore,
$$a(w)+b(w)+2c(w) = a(\swap_{(0,\tau)}(w))+b(\swap_{(0,\tau)}(w))+ 2c(\swap_{(0,\tau)}(w)),$$
and it follows that $L(w) = L(\swap_{(0,\tau)}(w)).$

\vskip.1in
\noindent
{\sf The case of $(\beta,\tau) \in \Swap^k_\twoswap(w)$.}

Of course $a(w) = a(\swap_{(\beta,\tau)}(w))$ because $w_{i_w(j), j} = w_{i_{\swap_{(\beta,\tau)}(w)}(j), j}$ 
for all column indices $j \in [n]$.
Using the abc Lemma ~\ref{abc_lemma}, to complete the proof, it suffices to show \sloppy{$b(w)+2c(w) = b(\swap_{(\beta,\tau)}(w))+2c(\swap_{(\beta,\tau)}(w)).$}

Observe that if $\{j,j'\} \cap \{\tau,\beta\} = \emptyset,$ then 
\begin{align}
b_{j,j'}(w) &= b_{j,j'}(\swap_{(\beta,\tau)}(w)) \label{swap_preserve_b}\\
c_{j,j'}(w) &= c_{j,j'}(\swap_{(\beta,\tau)}(w)),\label{swap_preserve_c}
\end{align}
because both $i_w(j)= i_{\swap_{(\beta,\tau)}(w)}(j)$ and $i_w(j') = i_{\swap_{(\beta,\tau)}(w)}(j')$.

\begin{slem}
\label{swapping_sublemma}
Let $w \in C_n, (\beta,\tau) \in \Swap^k_\twoswap(w).$ For any column index $j \in [n] - \left\{ \beta,\tau \right\},$ define
$$L_{\beta,\tau,j}(w) = \sum_{\substack{(x,y) \in [n]^2 \\ j \in \{x,y\} \\ \{x,y\} \cap \{\beta,\tau\} \neq \emptyset}}^{} b_{x,y}(w) + 2c_{x,y}(w).$$ Then,
\begin{equation}
\label{column_sum_invariance}
	L_{\beta,\tau,j}(w) = L_{\beta,\tau,j}(\swap_{(\beta,\tau)}(w)).
\end{equation}
\end{slem}
\sssec*{Proof of Lemma ~\ref{swapping_lemma} assuming Sublemma ~\ref{swapping_sublemma}}
As noted above, it suffices to show $b(w)+2c(w) = b(\swap_{(\beta,\tau)}(w))+2c(\swap_{(\beta,\tau)}(w)).$ By Sublemma ~\ref{swapping_sublemma},

\begin{align*}
	b(w)+2c(w)&=\sum_{j=1}^{n}\sum_{j'=1}^{n}(b_{j,j'}(w)+2c_{j,j'}(w))\\
	&=\sum_{\{j,j'\} \cap \{\tau,\beta\} = \emptyset}^{} (b_{j,j'}(w) +2c_{j,j'}(w))+ \sum_{\substack{1\leq j \leq n\\ j \notin \{\beta,\tau\}}}L_{\beta,\tau,j}(w)\\
	&= \sum_{\{j,j'\} \cap \{\tau,\beta\} = \emptyset}^{} (b_{j,j'}(\swap_{(\beta,\tau)}(w)) +2c_{j,j'}(\swap_{(\beta,\tau)}(w)))
	\\
	&\hspace{3 mm}+ \sum_{\substack{1\leq j \leq n\\ j \notin \{\beta,\tau\}}}L_{\beta,\tau,j}(\swap_{(\beta,\tau)}(w))\\
	&=\sum_{j=1}^{n}\sum_{j'=1}^{n}(b_{j,j'}(\swap_{(\beta,\tau)}(w))+2c_{j,j'}(\swap_{(\beta,\tau)}(w)))\\
	&= b(\swap_{(\beta,\tau)}(w))+2c(\swap_{(\beta,\tau)}(w)).\qedhere
\end{align*}
\end{proof}

\begin{proof}[Proof of Sublemma ~\ref{swapping_sublemma}]
We prove this in cases, depending on whether $i(j)$ is above, below, or row-between $i(\beta)$ and $i(\tau)$. The case $i(j)$ is row-between $i(\beta)$ and $i(\tau)$ is proven by considering the three further sub-cases that $j$ is right, left, or column-between $\beta$ and $\tau.$
\vskip.1in
\noindent
{\sf Case 1: Row $i(j)$ is above both $i(\beta)$ and $i(\tau).$}

If $i_w(j)$ is above both $i_w(\beta)$ and $i_w(\tau)$, hence $i_{\swap_{(b,t)}(w)}(j)$ is above both $i_{\swap_{(b,t)}(w)}(\beta)$ and $i_{\swap_{(b,t)}(w)}(\tau)$, we have $b_{x,y}(w)=b_{x,y}(\swap_{(b,t)}(w))$ and $c_{x,y}(w)=c_{x,y}(\swap_{(b,t)}(w))$ for all $(x,y) \in [n]^2$ with $j \in \{x,y\}, \{x,y\} \cap \{\beta,\tau\} = \emptyset.$ Hence, \eqref{column_sum_invariance} holds in this case.

\vskip.1in
\noindent
{\sf Case 2: Row $i(j)$ is below both $i(\beta)$ and $i(\tau).$}

This case is analogous to case 1.

\vskip.1in
\noindent
{\sf Case 3: Row $i(j)$ is row-between $i(\beta)$ and $i(\tau).$}
\newline\noindent
{\sf Sub-case 3a: Column $j$ is right of $\beta$ and $\tau.$}

If $j$ is right of $\beta$ and right of $\tau,$ then
 $$L_{\beta,\tau,j}(w) = b_{\beta,j}(w)+b_{\tau,j}(w)+2c_{\beta,j}(w)+2c_{\tau,j}(w),$$ as the other terms are 0. Now, for definiteness, assume $i(\beta)<i(\tau),$ and $j \not \equiv \beta$, as the other case is analogous. Then, 
$$b_{\beta,j}(w) = 0,b_{\tau,j}(w)=1,b_{\beta,j}(\swap_{(\beta,\tau)}(w)) = 1,b_{\tau,j}(\swap_{(\beta,\tau)}(w))=0.$$
Therefore,
$$b_{\beta,j}(w)+b_{\tau,j}(w)=b_{\beta,j}(\swap_{(\beta,\tau)}(w))+b_{\tau,j}(\swap_{(\beta,\tau)}(w)).$$
Analogously,
$$c_{\beta,j}(w)+c_{\tau,j}(w)=c_{\beta,j}(\swap_{(\beta,\tau)}(w))+c_{\tau,j}(\swap_{(\beta,\tau)}(w)),$$
and so
$L_{\beta,\tau,j}(w) = L_{\beta,\tau,j}(\swap_{(\beta,\tau)}(w))$ when $j$ is right of $\beta$ and right of $\tau$.
\newline\noindent
{\sf Sub-case 3b: Column $j$ is left of $\beta$ and $\tau.$}

This sub-case is analogous to sub-case 3a.
\newline\noindent
{\sf Sub-case 3c: Column $j$ is column-between $\beta$ and $\tau.$}

In this sub-case, if $j \equiv \beta$, we would obtain $L_{\beta,\tau,j} = 0$. So, assuming $j \not \equiv \beta$, we have
$$L_{\beta,\tau,j}(w) = b_{\beta,j}(w)+b_{j,\tau}(w)+2c_{\beta,j}(w)+2c_{j,\tau}(w),$$ and
\begin{align}
\label{l_equal_2}
L_{\beta,\tau,j}(w)=2=L_{\beta,\tau,j}(\swap_{(\beta,\tau)}(w)),
\end{align}
using the fact that $w_{i(j),j} \neq w_{i(\tau),\tau}.$ We check this in the case that $j \not \equiv \beta,i(\beta) > i(\tau),w_{i(\beta),\beta} = w_{i(\tau),\tau} = 1,$ and $w_{i(j),j} = -1$. Here,
$$b_{\beta,j}(w) = 1,b_{j,\tau}(w) = 1,c_{\beta,j}(\swap_{(b,t)}(w)) = 1 $$
and all other terms are $0$, yielding ~\eqref{l_equal_2}. The other cases are similar.
\end{proof}

\ssec{The Swapping Involution}
\label{the_swapping_involution}

We have now shown that $(b,t) \in \Swap^k(w)$ interacts nicely with $L$ and $\sgn.$ These properties will allow us to cancel out many terms in $S_{n,I}(X).$ However, in order to do so, we will define several involutions on $C_n^{I^c}$.
In order to define these involutions, we first need to define the swapping ordering.
This ordering is defined on pairs of columns so that one pair is less than another if the corresponding pair of rows is lexicographically higher in the matrix.
Then, for $\alpha \in \left\{ l,m,p \right\}$
we define involutions $\involution^k_\alpha$ given by applying the swap in $\Swap^k_\alpha$ minimum under the swapping ordering.

\begin{defn}
\label{swapping_ordering}
Let $w \in C_n.$ Define a total ordering, called the {\it swapping ordering}, on $\Swap^k(w)$ by $(p,q) < (b,t)$ if $(i(p), i(q))$ comes lexicographically before
$(i(b),i(t))$. More formally $(p,q) < (b,t)$ if either of the following holds:
\begin{itemize}
	\item $\min(i(p),i(q)) < \min(i(b),i(t))$ or
	\item $\min(i(p),i(q)) = \min(i(b),i(t)) \text{ and } \max(i(p),i(q)) < \max(i(b),i(t)).$
\end{itemize}
\end{defn}
\begin{defn}
\label{swapping_involution_defn}
For $\alpha \in \left\{ l,k,m \right\}$, define the map $\involution^k_\alpha:C_n^{I^c}\rightarrow C_n^{I^c},$ as follows. Let $w \in C_n^{I^c}.$ If $\Swap^k_\alpha(w)$ as defined in Definition ~\ref{general_swap_defn} is nonempty, let $(b,t) \in \Swap^k_\alpha(w)$ be the minimum element of $\Swap^k_\alpha(w)$ taken with respect to the swapping ordering defined in Definition ~\ref{swapping_ordering}. Then, $\involution^k_\alpha$ is given by
$$
w \mapsto \begin{cases}
	w &\text{ if } \Swap^k_\alpha(w) = \emptyset \\
	\swap_{(b,t)}(w) &\text{ otherwise.}
\end{cases}$$
\end{defn}

We next show the above maps are sign-reversing involutions which preserve $L$.
\begin{lem}
\label{involution}
For each $\alpha \in \left\{ l,p,m \right\}$, the map $\involution^k_\alpha$ is an involution on $C_n^{I^c}$. If $\Swap^k_\alpha(w) \neq \emptyset$ then
\begin{align*}
	L(w)&=L(\involution^k_\alpha(w))   \\
	\sgn(w) &= -\sgn(\involution^k_\alpha(w)).
\end{align*}
\end{lem}
\begin{proof}

Let us prove the statement for the case of $\alpha = m$, since the proof of the case of $\alpha = p$ is analogous. The proof for the case of $\alpha =l$ follows since $\Swap^k_l(w) \subset \Swap^k_{m}(w).$

By the Swapping Lemma, $L(\involution^k_m(w)) = L(w)$ and $\sgn(w) = -\sgn(\involution^k_m(w))$. So, it suffices to show that $\involution^k_m$ as defined is indeed an involution. Since $w = \swap_{(b,t)}(\swap_{(b,t)}(w))$, it suffices to show that if $(b,t)$ is minimum (under the swapping ordering of Definition ~\ref{swapping_ordering})
	in $\Swap^k_m(w),$ and $(p,q)$ is minimum under the swapping ordering in $\Swap^k_m(\swap_{(b,t)}(w)),$ then $(b,t) = (p,q).$ (This statement is plausible as $(b,t) \in \Swap^k_m(\swap^k_{(b,t)}(w))$ by Lemma ~\ref{swap_after_swap}.) Since $\Swap^k(w) = \Swap^k_\signswap(w) \cup \Swap^k_\twoswap(w),$ it suffices to show that none of the four cases below can occur, as this implies that the minimum $(p,q) \in \Swap^k_m(\swap_{(b,t)}(w))$ is not less than the minimum $(b,t) \in \Swap^k_m(w)$.

\vskip.1in
\noindent
{\sf Case 1: The case $(b,t)\in \Swap^k_\signswap(w)$, and $(p,q) \in \Swap^k_\signswap(\swap_{(b,t)}(w))$ with $(p,q)<(b,t)$.}

For this case, $b = p = 0$ by assumption. For any $v \in C_n^{I^c},|\Swap^k_\signswap(v)| \leq 1,$ by Property (3) of Definition ~\ref{swap_defn}. Since $(0,t) \in \Swap^k_\signswap(\swap_{(0,t)}(w)),$ we have $\{(0,t)\} = \Swap^k_\signswap(\swap_{(0,t)}(w)),$ so $q = t$.

\vskip.1in
\noindent
{\sf Case 2: The case $(b,t) \in \Swap^k_\signswap(w)$, and $(p,q) \in \Swap^k_\twoswap(\swap_{(b,t)}(w))$ with $(p,q)<(b,t)$.}

This cannot happen since $k$-sign swaps are always less than $k$-two swaps in the swapping ordering. 
Indeed, $(b,t) \in \Swap^k_\signswap(w)$ so $b = 0.$ However, since $p > 0$ which implies $(b,t) < (p,q).$

\vskip.1in
\noindent
{\sf Case 3: The case $(b,t) \in \Swap^k_\twoswap(w)$, and $(p,q) \in \Swap^k_\signswap(\swap_{(b,t)}(w)),$ with $(p,q) < (b,t)$.}

By assumption, $(p,q) \in \Swap^k_\signswap(\swap_{(b,t)}(w))$ so $p = 0$. 
First, we rule out the possibility that $(b,t) \in \Swap^k_{\dminus}(w) \cup \Swap^k_{\dplus}(w)$. 
To do this, note that because both $b$ and $t$ are greater than $k$, the set of rows $r$ with $j(r)$ left of $k$ in
$w$ and $\swap_{(b,t)}(w)$ are the same.
This implies that 
$\Swap^k_\signswap(\swap_{(b,t)}(w)) = \Swap^k_\signswap(w)$.
However, $\Swap^k_\signswap(w) = \emptyset$ because $(b,t)$ is the minimum element of $\Swap^k(w)$ under the swapping ordering, and any $k$-sign swap 
would be lower in the swapping ordering than $(b,t)$.
Therefore, $\Swap^k_\signswap(\swap_{(b,t)}(w)) = \emptyset$, contradicting that $(p,q) \in \Swap^k_\signswap(\swap_{(b,t)}(w))$.

To conclude this case, it suffices to rule out the subcase that $(b,t) \in \Swap^k_\leftswap(w) \cup \Swap^k_\sminus(w)$.
First, we must have $q \in \left\{ b,t \right\}$ since there are no column indices $s$ such that $i_w(s)$ is row-between $i_w(t)$ and $i_w(b)$ and $s \leq k$. 

If $(b,t) \in \Swap^k_\sminus(w)$, we see that $q = b$ and it follows from the definitions of single-minus $k$-swaps and $k$-sign swaps that
$(0,b) \in \Swap^k_m(w)$.
This contradicts that $(b,t)$ is minimum under the swapping ordering.

Finally, if $(b,t) \in \Swap^k_\leftswap(w)$, there are two cases
depending on whether
$i_w(b)>i_w(t)$ or 
$i_w(b)<i_w(t)$.
If $i_w(b)>i_w(t)$, we see $q = b$. Since $(0,b) \in \Swap^k_\signswap(\swap_{(b,t)}(w))$ it follows that $(0,t) \in \Swap^k_\signswap(w)$.
But $(0,t) < (b,t)$, contradicting that $(b,t)$ is minimum in the swapping ordering.
The case that 
$i_w(b)<i_w(t)$ is analogous with the roles of $b$ and $t$ reversed.
This completes the proof of case $3$.

\vskip.1in
\noindent
{\sf Case 4: The case $(b,t) \in \Swap^k_\twoswap(w)$, and $(p,q) \in \Swap^k_\twoswap(\swap_{(b,t)}(w))$ with $(p,q)<(b,t)$.}

If $(b,t)\in \Swap^k_{\dminus}(w) \cup \Swap^k_\dplus(w) \cup \Swap^k_\leftswap(w)$, then such a situation is impossible since the set of rows and signs of the rows $r$ with $j(r)$ left of $k$ (respectively right of $k$) in $\swap_{(b,t)}(w)$ and in $w$ are the same. Hence $\Swap^k(\swap_{(b,t)}(w))=\Swap^k(w)$. 

So, the only remaining case is when $(b,t)\in \Swap^k_{\sminus}(w)$. We first deal with the case $i_w(b)<i_w(t),$ $(p,q) \in \Swap^k_\twoswap(\swap_{(b,t)}(w)),$ and $(p,q)<(b,t).$ The only possibility is that one of $p,q$ is equal to $t,$ and $(p,q) \in \Swap^k_\sminus(\swap_{(b,t)}(w)) \cup \Swap^k_\dminus(\swap_{(b,t)}(w))$. For concreteness, let us say $q = t,$ the other case is the same with $p$ and $q$ reversed. If $(p,q) \in \Swap^k_\sminus(\swap_{(b,t)}(w))$, then one of $(p,b)$ or $(b,p)$ lies in $\Swap^k_\leftswap(w),$ and is less than $(b,t),$ contradicting the assumption that $(b,t)$ was minimal in $\Swap^k_m(w).$ If instead, $(p,t) \in \Swap^k_\dminus(\swap_{(b,t)}(w))$, then one of $(b,p)$ or $(p,b)$ would have been in $\Swap^k_\sminus(w)$. Again, this contradicts the assumption that $(b,t)$ was the minimal element of $\Swap^k_m(w)$ under the swapping ordering.

An analogous argument goes through in the case $i_w(b) > i_w(t),$ with $b,t$ reversed.
\end{proof}

\section{Part 1 of Theorem ~\ref{induction_result}}
\label{case_1}
Using the theory of swapping, we are able to prove Part 1 of Theorem ~\ref{induction_result} in Corollary ~\ref{part_1_induction_result}.

\begin{note}
\label{i_m_note}
Recall $n \in \BN,I = \{i_1,\ldots,i_l\}_< \subset [n-1]_0,i_0 = 0,$ and $i_{l+1} = n.$
For Section ~\ref{case_1} and Section ~\ref{case_2}, reserve the use of column index $i_m$ for the column index of the same name defined in Theorem ~\ref{induction_result}. 
In the case $m = i_m = -1$ Parts 1 and 2 of Theorem ~\ref{induction_result} hold automatically.
Therefore, in Sections ~\ref{case_1} and ~\ref{case_2}, we will assume $I$ is chosen so that $m \geq 0$.
\end{note}
Next, we define a $k$-block in $w$, which informally is a maximal contiguous set of rows, all of which have entries right of $k$.
\begin{defn}
\label{block_defn}
Let $w \in B_n,$ $k \in [n]$ be a column index, and $s,r \in [n]$ row indices with $r<s.$ Denote
$B_{r,s} = \{t \mid r \leq t \leq s\}.$ Call $B_{r,s}$ a {\it $k$-block in $w$} if it satisfies the following conditions:
\begin{itemize}
	\item For all $t \in B_{r,s},$ we have $j_w(t)>k.$
	\item If $s \neq n,$ then $j_w(s+1) \leq k.$
	\item If $r \neq 1,$ then $j_w(r-1) \leq k.$
\end{itemize}

If $B_{r,s}$ is a $k$-block in $w$ with $r \not\equiv s,$ call $B_{r,s}$ an {\it even $k$-block in $w$}, since it contains an even number of rows. If instead $r \equiv s,$ call $B_{r,s}$ an {\it odd $k$-block in $w$}.
Notate the set of all $k$-blocks in $w$ (respectively odd $k$-blocks in $w$, even $k$-blocks in $w$) by $\block^k(w)$ (respectively $\oblock^k(w),\eblock^k(w)$).
\end{defn}

\begin{lem}
\label{case_1_odd_differences}
Assume there is some $i_m \in I \cup \{i_{l+1}\}$, as defined in Notation ~\ref{i_m_note}. Let $i_b \in I \cup \{i_{l+1}\}$ with $i_b < i_m$. For all $w \in (C_n^{I^c})_{i_b},$ then $\oblock^{i_m}(w) \neq \emptyset$.
\end{lem}
\begin{proof}
Say $\block^{i_m}(w) = \{B_{r_1,s_1},\ldots, B_{r_h,s_h}\}$. The total number of rows in all the $i_m$-blocks in $w$ is equal to the number of columns right of $i_m$. For any $p$ with $1 \leq p \leq h$, we have $|B_{r_p,s_p}| = r_p - s_p+1,$ and the number of columns right of $i_m$ is $n - i_m.$ So, $\sum_{p=1}^{h}(s_p - r_p+1) = n-i_m.$ Furthermore, since $n-i_m = \delta_m + \delta_{m+1} + \cdots + \delta_l \equiv \delta_m$ is odd, there exists a $p$ for which $s_p - r_p+1$ is odd, and so $\oblock^{i_m}(w) \neq \emptyset$.
\end{proof}

\begin{lem}
\label{even_differences}
Assume there is some $i_m \in I \cup \{i_{l+1}\}$, as defined in Notation ~\ref{i_m_note}. If $i_k \in I \cup \{i_{l+1}\} ,w\in (C_n^{I^c})_{i_k},$ and $B_{r,s} \in \oblock^{i_m}(w)$, with $i_k < i_m,$ then $\Swap^{i_m}_l(w) \neq \emptyset$ and for all $(b,t) \in \Swap^{i_m}_l(w)$, there exists $h \in [n]_0$ with $h < i_m$ and $\swap_{(b,t)}(w) \in (C_n^{I^c})_h.$
\end{lem}
\begin{proof}
Let $B_{r,s} \in \oblock^{i_m}(w)$. First, we claim $s \neq n$. Since $w \in (C_n^{I^c})_{i_k},$ either $i(i_k) = n$ or $i(i_k +1) = n$. This is equivalent to $j(n) = i_k$ or $j(n) = i_k+1$. Thus, $j(n) \leq i_m$, implying $s \neq n$, as $j(s) > i_m$.

If $r \neq 1,$ then, since $s < n$, we obtain $s +1 \in [n],$ so $(j(r-1),j(s+1)) \in \Swap^{i_m}_\leftswap(w),$ or $(j(s+1),j(r-1)) \in \Swap^{i_m}_\leftswap(w).$ If $r = 1,$ then $(0,j(s+1)) \in \Swap^{i_m}_\signswap(w)$. Therefore, $\Swap^{i_m}_l(w) \neq \emptyset$.

Next, we show the second statement. First, we check $(j(n),i_m)\notin \Swap^{i_m}_l(w).$ If instead $(j(n),i_m) \in \Swap^{i_m}_l(w),$ then $j(n) \equiv i_m \not\equiv n,$ contradicting the assumption that $j(n) \equiv n.$
Hence, if $(b,t)$ is so that $t = i_m$, then either $\swap_{(b,t)}(w) \in (C_n^{I^c})_{j_w(n)}$ or $\swap_{(b,t)}(w) \in (C_n^{I^c})_{j_w(n)-1}$ because $b \neq j(n)$. Hence, $j_{\swap_{(b,t)}(w)}(n) = j_w(n)< i_m.$ If $t < i_m,$ then either $j_{\swap_{(b,t)}(w)}(n) = j_w(n)< i_m$ or $j_{\swap_{(b,t)}(w)}(n) = t< i_m,$ so again $\swap_{(b,t)}(w) \in (C_n^{I^c})_h$ with $h < i_m.$
\end{proof}

\begin{lem}
\label{equivalence_even_differences}
Assume there is some $i_m \in I \cup \{i_{l+1}\}$, as defined in Notation ~\ref{i_m_note}. Let $w \in C_n^{I^c}$. For any $h$, we have $\oblock^h(w) = \emptyset \implies \Swap^h_l(w) = \emptyset$. Consequently, if $w \in (C_n^{I^c})_{i_k}$ with $i_k < i_m,$ then $\oblock^{i_m}(w) = \emptyset \iff \Swap^{i_m}_l(w) = \emptyset.$
\end{lem}
\begin{proof}
Assume $\oblock^{h}(w) = \emptyset$. First, suppose $(b,t) \in \Swap^{h}_\leftswap(w).$ Then, if $i(b) < i(t)$, we have $B_{i(b) +1, i(t) - 1} \in \oblock^{h}(w)$, while if $i(t) < i(b),$ we have $B_{i(t) +1, i(b) - 1} \in \oblock^{h}(w)$. Therefore, $\Swap^{h}_\leftswap(w) = \emptyset$. Next, if $(b,t) \in \Swap^{h}_\signswap(w),$ then $B_{1,i(t) -1} \in \oblock^{h}(w)$. Therefore, $\Swap^{h}_\signswap(w) = \emptyset.$ Hence, $\Swap^{h}_l(w) = \Swap^{h}_\leftswap(w) \cup \Swap^{h}_\signswap(w) = \emptyset$.

For the second statement, the forward implication holds by the preceding paragraph, while the reverse direction holds by Lemma ~\ref{even_differences}.
\end{proof}

\begin{cor}
\label{part_1_induction_result}
Part 1 of Theorem ~\ref{induction_result} holds.
\end{cor}
\begin{proof}
Let $w \in (C_n^{I^c})_{i_b},$ with $i_b < i_m.$ Using Lemma ~\ref{case_1_odd_differences}, $\oblock^{i_m}(w) \neq \emptyset.$ But then, by Lemma ~\ref{equivalence_even_differences}, $\Swap^{i_m}_l(w) \neq \emptyset$, and so by Lemma ~\ref{involution}, the map $\involution^{i_m}_l$ defines a sign reversing involution on all such elements that additionally preserves $L$. Observe that this involution from Lemma ~\ref{involution} restricts to an involution on the set $\{w \in C_n^{I^c} \mid w \in (C_n^{I^c})_{i_b} \text{ with }i_b < i_m\}$ by the last part of Lemma 
\ref{even_differences}. Therefore,
\begin{align*}
	\sum_{i_k<i_m}^{} & S_{n,I,i_k}(X) = \frac{1}{2} \sum_{i_k < i_m}^{}\left(\sum_{w \in (C_n^{I^c})_{i_k}}^{}(-1)^{l(w)}L(w) + \sum_{w \in (C_n^{I^c})_{i_k}}^{}(-1)^{l(w)}L(w)\right)\\
	&= \frac{1}{2}\sum_{i_k<i_m}^{} \left(\sum_{w \in (C_n^{I^c})_{i_k}}^{}(-1)^{l(w)}L(w) +\sum_{w \in (C_n^{I^c})_{i_k}}^{}(-1)^{l(\involution^k_l(w))}L(\involution^k_l(w))\right)\\
	&= \frac{1}{2}\sum_{i_k<i_m}^{} \left(\sum_{w \in (C_n^{I^c})_{i_k}}^{}(-1)^{l(w)}L(w) -\sum_{w \in (C_n^{I^c})_{i_k}}^{}(-1)^{l(w)}L(w)\right)\\
	 &= 0.\qedhere
\end{align*}
\end{proof}

\section{Part 3 of Theorem ~\ref{induction_result} Implies Part 2 of Theorem ~\ref{induction_result}}
\label{case_2}
In this section, we present a bijection between 
$(C^{I^c}_n)_{i_m}$ and \sloppy{$(C^{J^c}_n)_{i_m+1}$}, where $J=\{i_1,\ldots, i_{m-1},i_m+1,i_{m+1},\ldots,i_l\},$
given by flipping the sign of the entry in the bottom row. Using direct computations in Lemma ~\ref{case_2_L} we show that both $L$ and $\sgn$ change nicely under this bijection, which will allow us reduce the proof Part 2 of Theorem ~\ref{induction_result} to Part 3 of Theorem ~\ref{induction_result}. After this section, we prove Part 3 of Theorem ~\ref{induction_result}.

\begin{lem}
\label{case_2_L}
Assume there is some $i_m \in I$, as defined in Notation ~\ref{i_m_note},
and let $J=\{i_1,\ldots, i_{m-1},i_m+1,i_{m+1},\ldots,i_l\}$. Define maps 
\begin{align*}
	g:(C^{I^c}_n)_{i_m} &\longrightarrow (C^{J^c}_n)_{i_m+1} \\
	w &\longmapsto s_{n-1}s_{n-2} \cdots s_1 s_0 s_1 \cdots s_{n-1}w, \\
	f:(C^{J^c}_n)_{i_m+1} &\longrightarrow (C^{I^c}_n)_{i_m}\\
	w &\longmapsto s_{n-1}s_{n-2} \cdots s_1 s_0 s_1 \cdots s_{n-1}w,
\end{align*}
where $s_0,\ldots, s_{n-1}$ are the Coxeter generators of $B_n$. Then, $f,g$ are mutual inverses and for any $w \in (C^{I^c}_n)_{i_m},$
we have $\sgn(g(w)) = -\sgn(w)$ and $L(g(w))+i_m +1 = L(w).$
\end{lem}
\begin{proof}
First, note that $g$ indeed sends $(C^{I^c}_n)_{i_m}$ to $(C^{J^c}_n)_{i_m+1}$ since it sends the matrix $w$ to a signed permutation matrix which only differs from $w$ in that the entry in the bottom row is changed from a ``$-1$'' to a ``$1$''. The fact that $f$ and $g$ are mutual inverses is immediate from the definition. The fact that $\sgn(g(w)) = -\sgn(w)$ holds since $g$ is multiplication by an odd number of transpositions.

Finally, we check $L(g(w)) + i_m + 1 = L(w)$. Note that since $n - i_m$ is odd, we have $w_{n,i_m+1}=-1$ by Definition ~\ref{cnik_defn}.
We now use the abc Lemma ~\ref{abc_lemma} to compute $L.$ First, because $b_k(w) = b_k(g(w))$ for all $k \in [n],$ we have $b(w)= b(g(w)).$

To compute $a(g(w))$ and $c(g(w))$ in terms of $a(w),c(w),$ and $i_m$, we have two cases: $n$ is even or odd. We will only check the even case, since the odd case is similar. If $n$ is even, then $i_m+1$ is even, which means there are $\frac{i_m+1}{2}$ odd columns left of $i_m+1.$ This tells us that $c_{i_m+1}(w) = \frac{i_m+1}{2}$ and $c_{i_m+1}(g(w)) = 0.$ Additionally, $c_k(w) = c_k(g(w))$ for $k \neq i_m+1.$ Also, since $i_m+1$ is even, $a(w) = a(g(w)).$ This tells us that 
\begin{align*}
	L(w) - L(g(w)) &= a(w) - a(g(w)) + b(w) - b(g(w)) +2c(w) - 2c(g(w)) 
	\\
	&= 0 + 0 + 2c_{i_m+1}(w)
	\\
	&= 2\cdot \frac{i_m+1}{2} 
	\\
	&= i_m+1.
\end{align*}
Hence, $L(g(w))+i_m +1 = L(w).$
\end{proof}

\begin{lem}
\label{part_2_induction_result}
Assuming Part 3 of Theorem ~\ref{induction_result}, Part 2 holds.
\end{lem}
\begin{proof}
Let $J=\{i_1,\ldots, i_{m-1},i_m+1,i_{m+1},\ldots,i_l\}.$
Using the function $g$ as defined in Lemma ~\ref{case_2_L},
\begin{align*}
	\sum_{w \in (C_n^{I^c})_{i_m}}^{} (-1)^{l(w)}X^{L(w)}&=\sum_{w \in (C_n^{I^c})_{i_m}}^{} -(-1)^{l(g(w))}X^{L(g(w))+i_m+1} & \text{(by Lemma ~\ref{case_2_L})}\\
	&= \sum_{u\in (C_n^{J^c})_{i_m+1}}^{} -(-1)^{l(u)}X^{L(u)+i_m+1}& \text{(by Lemma ~\ref{case_2_L})}
	\\
	&=-X^{i_m+1}\sum_{u\in (C_n^{J^c})_{i_m+1}}^{} (-1)^{l(u)}X^{L(u)}
	\\
	&= -X^{i_m-1}X^{n-(i_m-1)} f_{n,J^{(m)}}(X)& \text{(by Theorem ~\ref{induction_result} Part 3)}
	\\
	&= -X^nf_{n,J^{(m)}}(X)
	\\
	&= -X^nf_{n,I^{(m+1)}}(X) & \text{(since $J^{(m)} = I^{(m+1)}$).}
\end{align*}
\end{proof}

\section{Part 3 of Theorem ~\ref{induction_result}}
\label{case_3}

In this section, we prove Part 3 of Theorem ~\ref{induction_result}. The idea of the proof is to
define a bijection between elements $w \in (C^{I^c}_n)_j$ with $\Swap^j_m(w) = \emptyset$ and elements of $v \in C^{(I^{(r)})^c}_n$ with $\Swap^{j-1}_p(v) = \emptyset$, where $j,r$ are defined in Notation ~\ref{fixed_j} below. This bijection is given in Algorithm ~\ref{algorithm}.

\begin{note}
\label{fixed_j}
Let $n \in \BN, I = \{i_1,\ldots, i_l\}_<, i_0 = 0,$ and $i_{l+1} = n$. 
Let $r \in [l+1]_0$ be so that $\delta_k \equiv 0$ whenever $r \leq k \leq l$.
Because this condition is vacuously satisfied when $r = l+1$, such an $r$ always exists.
{\bf For the remainder of the paper,} use $j$ to denote the fixed column index $i_r$. By Definition ~\ref{cnik_defn}, if $w \in (C_n^{I^c})_j$ we have $w_{i(j),j} = 1.$
\end{note}

\subsection{Initially Uncancelled Elements}

In order to define initially uncancelled elements in Definition ~\ref{initially_uncancelled_defn}, we first need some terminology.
To start, we define $(\varepsilon_1, \varepsilon_2)^k_t$ pairs, which are pairs of consecutive columns $(t,t+1)$ right of $k$ 
with $t \not \equiv n$ and nonzero entries
$\varepsilon_1, \varepsilon_2$.
More precisely:

\begin{defn}
\label{pairs}
Let $k \in [n]_0,w \in (C_n^{I^c})_k,$ and $t \in [n-1]$. For $\epsilon_1, \epsilon_2 \in \left\{ +1, -1 \right\}$ we say $w$ has an {\it $(\epsilon_1, \epsilon_2)^k_t$ pair} if
\begin{itemize}
	\item $(\epsilon_1, \epsilon_2) = (w_{i(t),t},w_{i(t+1),t+1})$,
	\item $t \not \equiv n$,
	\item $t > k$.
\end{itemize}
\end{defn}

\begin{rem}
\label{no_plus_minus_pairs}
Suppose $w \in (C_n^{I^c})_j$. Note that we cannot have any $t$ for which $w$ has a $(1,-1)^j_t$ pair, as for any column $t$ which is right of $j$ with $t \not \equiv n,$ we have $t \notin D(w)$. Hence, the only possibilities are $(1,1)^j_t,(-1,-1)^j_t$ and $(-1,1)^j_t$ pairs.
\end{rem}

We next define blocks of $1$s and blocks of $-1$s. Blocks of $1$s are subsets of blocks which are maximal sets of adjacent rows with all $1$s, while blocks of $-1$s are subsets of blocks which are maximal sets of adjacent rows with all $-1$s. More formally:

\begin{defn}
\label{defn:block_of_ones}
Let $w \in B_n.$ Let $s,r \in [n]$ be row indices with $r\leq s.$ Call the set of row indices $\plussetblock_{r,s}=\{t \mid r \leq t \leq s\}$ a {\it $k$-block of $1$s in $w$} if
\begin{itemize}
	\item For all $t \in B_{r,s},$ we have $w_{t, j(t)} = 1$ and $j(t)>k.$
	\item If $s \neq n,$ then $j(s+1) \leq k$ or $w_{s+1,j(s+1)} = -1.$
	\item If $r \neq 1,$ then $j(r-1) \leq k$ or $w_{r-1,j(r-1)} = -1.$ 
\end{itemize}
If $r \not \equiv s,$ then $\plussetblock_{r,s}$ is called an {\it even $k$-block of $1$s}. If $r \equiv s,$ $\plussetblock_{r,s}$ is called an {\it odd $k$-block of $1$s}.

Call the set of row indices $\minussetblock_{r,s}=\{t \mid r \leq t \leq s\}$ a {\it $k$-block of $-1$s in $w$} if
\begin{itemize}
	\item For all $t \in B_{r,s},$ we have $w_{t, j(t)} = -1$ and $j(t)>k.$
	\item If $s \neq n,$ then $j(s+1) \leq k$ or $w_{s+1,j(s+1)} = 1.$
	\item If $r \neq 1,$ then $j(r-1) \leq k$ or $w_{r-1,j(r-1)} = 1.$ 
\end{itemize}
If $r \not \equiv s,$ then $\minussetblock_{r,s}$ is called an {\it even $k$-block of $-1$s in $w$}, and, if $r \equiv s,$ then $\minussetblock_{r,s}$ is called an {\it odd $k$-block of $-1$s in $w$}.

Notate the set of all $k$-blocks of $1$s in $w$ (respectively even $k$-blocks of $1$s in $w$, odd $k$-blocks of $1$s in $w,$ $k$-blocks of $-1$s in $w$, even $k$-blocks of $-1$s in $w$, odd $k$-blocks of $-1$s in $w$) by $\pblock^k(w)$ (respectively $\epblock^k(w),$ $\opblock^k(w),\mblock^k(w), \emblock^k(w),\omblock^k(w)$).
\end{defn}

We next define initially uncancelled elements.
The seven conditions defining an initially uncancelled are quite opaque, and we recommend the reader does not try to memorize them, but instead
only refer to them as needed.
The reason we make this definition is that these seven conditions turn out to precisely characterize those elements in $(C_n^{I^c})_j$ with no general-minus $j$-swaps,
as is shown in Theorem ~\ref{equivalence_uncancelled} below.

\begin{defn}
\label{initially_uncancelled_defn}
Let $I = \{i_1,\ldots, i_l\}$ and fix $k \in I \cup \{n\}$ such that $k\equiv i$ for all $i \in I \cup \{n\}$ with $i>k.$
An element $w \in (C_n^{I^c})_k$ is {\it $k$-initially uncancelled} if its signed permutation matrix satisfies the following properties.
\begin{enumerate}
	\item Either $j(n-1)$ is left of $k$ or $w_{n-1,j(n-1)} = -1.$
	\item $\oblock^k(w) = \emptyset$.
	\item Suppose $t$ is a row index so that $w_{t,j(t)}=1,$ $j(t)$ is right of $k,$ and there is no $s \in [n]$ with $s \geq t$ so that $\plussetblock_{1,s} \in \opblock^k(w)$. Then, there is some $\plussetblock_{r,s} \in \epblock^k(w)$ with $r \leq t \leq s.$
	\item If $n$ is even, there is no $t \in [n-1]$ for which $w$ has a $(-1,1)^k_t$ pair. If $n$ is odd, there is at most one $t \in [n-1]$ for which $w$ has a $(-1,1)^k_t$ pair, and if $w$ has such a pair, then $i(t+1) = 1.$ Furthermore, $\plussetblock_{1,i(t)-1} \in \pblock^k(w).$
	\item If $t$ is so that $w$ has a $(-1,-1)^k_t$ pair, then either $i(t+1) +1 = i(t)$ or $B_{i(t+1)+1,i(t)-1} \in \epblock^k(w).$
	\item Suppose $\minussetblock_{p,q} \in \mblock^k(w)$ and $B_{r,s} \in \block^k(w)$ with $r>1$ and $\minussetblock_{p,q} \subsetneq B_{r,s}$. Then $\minussetblock_{p,q} \in \omblock^k(w)$ if and only if $p=r$ or $q = s.$
	\item If $r >1$ and $s$ are row indices with $B_{r,s} \in \block^k(w)$, and $q$ is a row index with $\plussetblock_{r,q} \in \pblock^k(w)$ (respectively $\plussetblock_{p,s} \in \pblock^k(w)$), then $q = s$ (respectively $p = r$).
	\end{enumerate}
We notate $(U_n^{I^c})_k= \{w \in (C_n^{I^c})_k \mid w \text{ is }k\text{-initially uncancelled}$\}.
\end{defn}

\begin{eg}
\label{initial_example}
Take
$$w=\begin{pmatrix}
	 &  &  &  & &  &  1&  &  & \\
	 &  &  &  &  &  &  &  1&  & \\
	 &  &  &  & 1 &  &  &  &  & \\
	 &  &  &  &  & 1 &  &  &  & \\
	 1&  &  &  &  &  &  &  &  & \\
	 &  &  & -1 &  &  &  &  &  & \\
	 &  &  &  &  &  &  &  & 1 & \\
	 &  &  &  &  &  &  &  & & 1\\
	 &  &  -1&  &  &  &  &  &  &\\
	 &  1&  &  &  &  &  &  &  & \\
\end{pmatrix}$$
Here, $w\in (U_{10}^{\{2,6\}^c})_2$. To see this, note that $D(w) = \{2,6\}.$ Property (1) holds because $w_{9,3} = -1.$ Property (2) holds because $B_{1,4} \in \eblock^2(w)$ and $B_{6,9} \in \eblock^2(w)$.  Property (3) holds because $\plussetblock_{7,8} \in \epblock^2(w)$ and rows $1,2,3,4$ also determine $\plussetblock_{1,4} \in \epblock^2(w).$ Property (4) is satisfied because there are no $t$ for which $w$ has a $(-1,1)^2_t$ pair.
Property (5) holds because columns $3,4$ are to the right of $2,$ and the only rows between $i(3)=9$ and $i(4) = 6$ are rows $7$ and $8$ which both contain $1$s right of $2.$ Property (6) holds because the only times $\minussetblock_{p,q} \subset B_{r,s},\minussetblock_{p,q} \neq B_{r,s}$ are when $(r,s) = (6,9)$ and either $p = q = 6$ or $p = q = 9,$ but in both cases $\minussetblock_{p,q} \in \omblock^2(w)$. Finally, Property (7) is satisfied because in $B_{1,4}$ there are no rows with a $-1,$ and in $B_{6,9},$ the top and bottom rows contain $-1$s.
\end{eg}
\ssec{Equivalence of Initially Uncancelled and an empty set of General-Minus $j$-Swaps}
We show in Theorem ~\ref{equivalence_uncancelled} that $w \in (C_n^{I^c})_j$ satisfies the seven properties of $j$-initially uncancelled matrices from Definition ~\ref{initially_uncancelled_defn} if and only if $\Swap_m^j(w) = \emptyset$. The if direction is much harder, since we must verify each of the seven properties of $j$-initially uncancelled matrices, one at a time, assuming $\Swap^j_m(w)=\emptyset$.
\sssec{Property (1)}

\begin{lem}
\label{n-1_row}
Let $w\in (C_n^{I^c})_j$. We cannot have both that column $j(n-1)$ is right of $j$ and $w_{n-1,j(n-1)}=1.$
\end{lem}
\begin{proof}
Suppose column $j(n-1)$ is right of $j$ and $w_{n-1,j(n-1)}=1$. Since $j(n-1) \not \equiv n,$ we have $j(n-1) \notin D(w)$. Hence, the only possibility is that 
$w_{n,j(n-1)+1} = 1,$
 contradicting the assumption that $j(n-1) > j$ and $w_{n,j} = 1$.
\end{proof}

The next lemma is crucial, as while Lemma ~\ref{involution} allows us to cancel out elements of $C_n^{I^c}$ with $\involution^j_m \neq \emptyset$ from $S_{n,I}(X)$, Lemma ~\ref{restrict_swaps} allows us to cancel out the contributions of elements from $(C_n^{I^c})_j$ with $\involution^j_m \neq \emptyset$ from $S_{n,I,j}.$ 

\begin{lem}
\label{restrict_swaps}
The involution $\involution^j_m:C_n^{I^c}\rightarrow C_n^{I^c}$ defined in \linebreak Lemma ~\ref{involution} restricts to an involution $\involution^j_m:(C_n^{I^c})_j\rightarrow (C_n^{I^c})_j.$
\end{lem}
\begin{proof}
Since we have shown that $\involution^j_m$ defines an involution on $C_n^{I^c}$ in Lemma ~\ref{involution}, it suffices to show that we have a containment $\involution^j_m(C_n^{I^c})_j \subset (C_n^{I^c})_j.$ That is, we will show  that if $w \in C_n^{I^c}$ and $(b,t)$ is the minimum element of $\Swap^j_m(w)$ under the swapping ordering, so that $\involution^j_m(w)= \swap_{(b,t)}(w),$ then $b \neq j$ and $t \neq j$. 

We show that $b \neq j, t \neq j$ by showing it for each of the different types of swaps.
First, suppose $(b,t) \in \Swap^j_l(w)$. Since $b < t \leq j$, we only need rule out the possibility $t = j.$ If $t = j$, 
$\swap_{(b,t)}(w)$ satisfies the hypotheses of Lemma ~\ref{even_differences}, so $\involution^{j}_l(\swap_{(b,t)}(w)) \notin (C_n^{I^c})_j.$ Hence, $\involution^{j}_l(\swap_{(b,t)}(w)) \neq w,$ contradicting Lemma ~\ref{involution}. Next, if \sloppy{$(b,t) \in \Swap^j_\dminus(w) \cup \Swap^j_\dplus(w)$}, then $j < b < t.$ So, in this case, $\swap_{(b,t)}(w) \in (C_n^{I^c})_j$. The final possibility is that $(j,t) =(b,t) \in \Swap^j_\sminus(w)$. Because $(j,t)$ is a $j$-swap, we must have that $j(n-1)$ is right of $j$ and $w_{n-1,j(n-1)} = 1,$ contradicting Lemma ~\ref{n-1_row}. Therefore, column $j$ will never be involved in the least swap of $w$, and so $\involution^j_m:C_n^{I^c}\rightarrow C_n^{I^c}$ restricts to a map $\involution^j_m:(C_n^{I^c})_j\rightarrow (C_n^{I^c})_j.$
\end{proof}

\sssec{Property (2)}

This was already shown in Lemma ~\ref{even_differences}.

\sssec{Property (3)}
\begin{lem}
\label{even_ones}
Any matrix $w\in (C_n^{I^c})_j$ with $\Swap^j_m(w) = \emptyset$ has 
$\opblock^j(w)= \emptyset$ or 
$\opblock^j(w) = \{\plussetblock_{1,s}\}$ for some $s \in [n].$
\end{lem}

\begin{proof}
	Suppose $\plussetblock_{r,s} \in \opblock^j(w),$ with $r \neq 1.$ Then, observe $(r-1,s+1) \in \Swap^j_\leftswap(w) \cup \Swap^j_\dminus(w) \cup \Swap^j_\sminus(w).$ Indeed, this holds because for for $\alpha \in \left\{ r-1, s+1 \right\}$ we must either have $j(\alpha)$ is left of $j$ or $w_{\alpha, j(\alpha)} = -1$, by definition of odd $j$-block of $1$s.
\end{proof}

\sssec{Some Intermediate Lemmas}
\begin{lem}
\label{no_odd_minus}
If $w\in (C_n^{I^c})_j$ has $\Swap^j_m(w) = \emptyset$ and $B_{r,s}\in \eblock^j(w)$, then there cannot be $\minussetblock_{p,q} \in \omblock^j(w)$ with $r<p<q<s.$
\end{lem}
\begin{proof}
If $\minussetblock_{p,q} \in \omblock^j(w)$ with $r<p<q<s$, then $(p-1,q+1) \in \Swap^j_{\dplus}(w)$.
\end{proof}
We next define statistics $s_e^k, s_o^k, u_e^k$ and $u_o^k$. The statistic $s_e^k$ (respectively $s_o^k$)
is the set even (respectively odd) indexed columns right of $k$ containing a $-1$.
The statistic $u_e^k$ (respectively $u_o^k$) is the set of columns $h$ right of $k$ containing a $-1$ with an even (respectively odd)
number of rows $t$ left of $k$ such that $i(t)$ is above $i(h)$.

\begin{defn}
\label{us_defn}
Let $w\in B_n,k \in [n]$ and define
$$s_e^k(w) = \{h \in [n] \mid  h>k,h \equiv 0,w_{i(h),h}=-1\},$$
$$s_o^k(w) = \{h \in [n] \mid  h>k,h \equiv 1,w_{i(h),h}=-1\}.$$
Also, for each column index $h \in [n]$ with $h>k$ and $w_{i(h),h} = -1,$ let 
$$up_h^k(w) = |\{t \in [k-1] \mid i(t)<i(h)\}|.$$
Then define 
$$u_e^k(w) = \{h \in [n] \mid  h>k,w_{i(h),h}=-1,|up_h^k(w)|\equiv 0\}$$
$$u_o^k(w) = \{h \in [n] \mid  h>k,w_{i(h),h}=-1,|up_h^k(w)|\equiv 1\}.$$
For $a,b \in \{o,e\}$ let $u_a^j(w)s_b^j(w) = u_a^j(w)\cap s_b^j(w).$

For $B_{r,s} \in \eblock^k(w),$ define $\minus_{B_{r,s}} = \{t\in B_{r,s} \mid w_{t,j(t)}=-1\}.$
\end{defn}

\begin{rem}
\label{us_and_pairs}
The subscript $e$ stands for even and the subscript $o$ stands for odd. These statistics are crucial for seeing how $L$ changes when we apply Algorithm \ref{algorithm}.
\end{rem}

\begin{lem}
\label{even_minus_ones}
Let $w\in (C_n^{I^c})_j$ with $\Swap^j_m(w) = \emptyset$. If $B_{r,s} \in \eblock^j(w)$ with $r \neq 1,$ then $|\minus_{B_{r,s}} \cap s_e^j(w)| =|\minus_{B_{r,s}}\cap s_o^j(w)|.$ If $r = 1,$ then
$$\left| |\minus_{B_{r,s}} \cap s_e^j(w)|-|\minus_{B_{r,s}}\cap s_o^j(w)| \right|=\begin{cases}
	 1, &\text{ if } |\minus_{B_{r,s}}| \equiv 1,\\
	0,  &\text{ if } |\minus_{B_{r,s}}| \equiv 0.
\end{cases}$$
\end{lem}
\begin{proof}
First, consider the case that $B_{r,s} \in \eblock^j(w)$ with $r \neq 1.$ Whenever some $t \in B_{r,s}$ satisfies $w_{t,j(t)} = 1,$ we know by Lemma ~\ref{even_ones} that there is some $\plussetblock_{a,b} \in \epblock^j(w)$ with $a \leq t \leq b.$ Write $\minus_{B_{r,s}} = \{k_1,\ldots, k_p\}_<.$ First, observe that by Lemma ~\ref{even_differences}, $p$ must be even, since $|B_{r,s}| \equiv 0$, and $|\{t \in B_{r,s} \mid w_{t,j(t)} = 1\}|\equiv 0,$ as $\{t \in B_{r,s} \mid w_{t,j(t)} = 1\}$ is a disjoint union of sets of even cardinality. So $|B_{r,s}| - |\{t \in B_{r,s} \mid w_{t,j(t)} = 1\}| =|\minus_{B_{r,s}}|$ is also even. Next, note that for $1 \leq i \leq p-1,$ either $k_i+1 = k_{i+1},$ or else $\plussetblock_{k_i+1,k_{i+1}-1} \in \epblock^j(w).$ This tells us that $k_i \not \equiv k_{i+1}.$ Therefore, for any $B_{r,s} \in \eblock^j(w)$ such that $r \neq 1,$ we must have 
$|\minus_{B_{r,s}} \cap s_e^j(w)| =|\minus_{B_{r,s}}\cap s_o^j(w)|.$ 

In the case that $r = 1,$ the argument from the previous paragraph shows that the parities of the row indices in $\minus_{B_{r,s}}$ must alternate. However, if $r=1,$ then $|\minus_{B_{r,s}}|$ may be either even or odd. Hence, $\left| |\minus_{B_{r,s}} \cap s_e^j(w)|-|\minus_{B_{r,s}}\cap s_o^j(w)| \right|$ will either be $0$ or $1$ depending on whether $|\minus_{B_{r,s}}|$ is even or odd, as claimed.
\end{proof}

\begin{lem}
\label{pair_setup}
Any $w\in (C_n^{I^c})_j$ with $\Swap^j_m(w) = \emptyset$ can have at most one $(-1,1)^j_t$ pair. Furthermore $w$ has a $(-1,1)^j_t$ pair if and only if there is some $q \in [n]$ for which $\plussetblock_{1,q} \in \opblock^j(w)$.
\end{lem}

\begin{proof}

For definiteness, let us take the case that $n \equiv 0.$ The case $n \equiv 1$ is exactly analogous, simply by reversing the roles of even and odd.

Since $n$ is even, if $w$ has a $(-1,1)^j_m$ pair, then $m$ is odd, because by Definition ~\ref{pairs} $m \not \equiv n.$ So, if there were any $(-1,1)^j_m$ pairs, we must have $u_o^j(w)>u_e^j(w),$ by Remark ~\ref{no_plus_minus_pairs}.

By Lemma ~\ref{even_minus_ones} we know that for any $B_{r,s} \in \eblock^j(w)$ with $r \neq 1,$ \sloppy{$|\minus_{B_{r,s}} \cap s_e^j(w)| =|\minus_{B_{r,s}}\cap s_o^j(w)|.$} By Remark ~\ref{no_plus_minus_pairs}, the number $t \in [n-1]$ for which $w$ has a $(-1,1)^j_t$ pair is exactly the sum of the differences
\begin{align*}
\sum_{B_{r,s} \in \block^j(w)}^{}\left(|\minus_{B_{r,s}} \cap s_o^j(w)|-|\minus_{B_{r,s}}\cap s_e^j(w)|\right). 
\end{align*}
By Lemma ~\ref{even_minus_ones}, this is $0$ if there are no $s$ for which $B_{1,s} \in \eblock^j(w),$ and is $|\minus_{B_{1,s}}| \bmod 2$
if there is such an $s$. This shows there is at most one $t$ for which $w$ has a $(-1,1)^j_t$ pair.

To see the second statement, note that $|\minus_{B_{1,s}}| \equiv 1$ if and only if $\opblock^j(w) \neq \emptyset$. By Lemma ~\ref{even_ones}, the only way there is some $\plussetblock_{p,q} \in \opblock^j(w)$ with $\plussetblock_{p,q} \subset\minus_{B_{1,s}}$ is if $p = 1.$ Hence, $w$ has a single $(-1,1)^j_t$ pair if and only if $\plussetblock_{p,q} \in \opblock^j(w)$ with $p = 1.$
\end{proof}

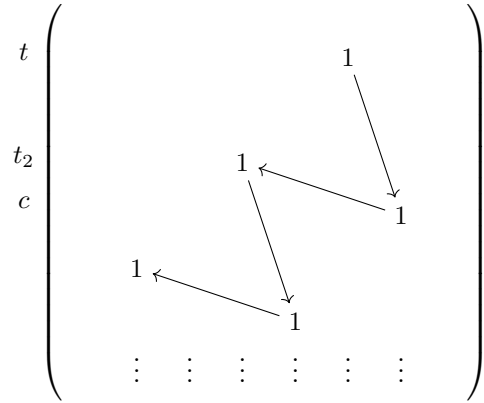
\begin{figure}[h!]
\[
	\begin{tikzpicture}[baseline=-0.5ex]
\matrix(B)[matrix of math nodes,column sep={20pt,between origins},row
    sep={20pt,between origins},left delimiter={(},right delimiter={)}] (m)
    {
    &  &  &  &  & & &  \\
	 &  &  &  &  & |[name=1]|1 &  &  \\
	 |[name=left]|&  &  &  &  &  &  &  \\
	 &  &  & |[name=3]|1 &  &  &  &  \\
	 &  &  &  & &  & |[name=2]| 1 &  \\
	 & |[name=5]|1 &  &  &  &  &  & \\
	 &  &  &  & |[name=4]|1 &  &  &  \\
	 & \vdots & \vdots& \vdots &  \vdots& \vdots &\vdots  & \\
  };
  \draw[->]
  (1) to (2);
  
  \draw[->]
  (2) to (3);
  
  \draw[->]
  (3) to (4);
  
  \draw[->]
  (4) to (5);
   \node at (-3.2,2) {$t$};
   \node at (-3.2,0) {$c$};
   \node at (-3.2,.6) {$t_2$};
  \end{tikzpicture}
  \]
\caption{A schematic representation of the infinite descent in the proof of Lemma ~\ref{adjacent_ones}.}
\label{infinite_descent_matrix}
\end{figure}

The next lemma uses the technique of infinite descent. See Figure ~\ref{infinite_descent_matrix} for a pictorial representation of the idea of proof.
\begin{lem}
\label{adjacent_ones}
Let $w\in (C_n^{I^c})_j$ with $\Swap^j_m(w) = \emptyset$. If $B_{r,s} \in \block^j(w)$ and $t\in B_{r,s}$ so that $t \not \equiv n$ and $w_{t,j(t)} = 1,$ then $w_{t+1,j(t)+1} = 1.$
\end{lem}

\begin{proof}
We know that column $j(t)+1$ must have a $1$ as $t \notin D(w)$ because $t \not \equiv n.$ So, to complete the proof, we only need to show that $j(t+1) = j(t)+1.$ Suppose instead $j(t+1) \neq j(t)+1.$
We will show that we will get stuck in an infinite descent. Let $c=i(j(t)+1).$ Then, $w_{c,j(c)} = 1,$ and either $w_{c-1,j(c-1)} = 1,$ or there is some $\plussetblock_{b_1,t_1} \in \epblock^j(w)$, with $c \in \plussetblock_{b_1,t_1},$ by Lemma ~\ref{even_ones}. So, either $w_{c-1,j(c-1)} = 1,$ and $j(c-1)$ is right of $j$ or else $c= b_1$.

Suppose we are in the former case, that $w_{c-1,j(c-1)} = 1,$ and $j(c-1)$ is right of $j.$ Define $t_2 = c-1.$ We know $j(t_2) \neq j(c)-1,$ since $j(t) = j(c)-1,$ and by assumption, $j(t)+1 \neq j(t+1).$ Therefore, since $t_2 \not\equiv n,$ we are again in the same situation as before, except $t$ is replaced by $t_2.$

The latter case that $c \in \plussetblock_{c,t_1}$ is handled similarly, by applying the argument of the preceding paragraph to $w_{t_1,j(t_1)}$ instead of $w_{c-1,j(c-1)}.$

In both cases, we can repeat the above argument indefinitely, and we will be forced to move down the matrix indefinitely. This implies $n$ is unbounded, which is a contradiction.
\end{proof}
\sssec{Property (4)}
\begin{lem}
\label{no_even_pairs}
Let $w\in (C_n^{I^c})_j$ with $\Swap^j_m(w) = \emptyset$. If $n$ is even, there is no $t \in [n-1]$ so that $w$ has a $(-1,1)^j_t$ pair.
\end{lem}
\begin{proof}

Assume $w$ has some $(-1,1)^j_t$ pair, and let $B_{1,s} \in \block^j(w)$. We saw in Lemma ~\ref{pair_setup} that there is some $\plussetblock_{1,q} \in \opblock^j(w)$ with $q \geq s$. However, we must then have $w_{q+1,j(q+1)} = -1$ and $j(q+1)$ is right of $j.$ This is because if $j(q+1)$ were left of $j,$ we would have $B_{1,s} = \plussetblock_{1,q}$ which would imply $B_{1,s} \in \oblock^j(w)$, which we cannot have by Lemma ~\ref{even_differences}

Therefore, $w_{q+1,j(q+1)} = -1.$ We know $q \not\equiv n \equiv 0,$ so $q \equiv 1 \not\equiv n.$ But then, by Lemma ~\ref{adjacent_ones}, if follows that $w_{q+1,i(q)+1} = 1,$ contradicting $w_{q+1,j(q+1)} = -1,$ since a signed permutation matrix can only have one nonzero entry in each row. Therefore, $w$ has no $(-1,1)^j_t$ pairs.
\end{proof}

\begin{lem}
\label{one_odd_pair}
Let $w\in (C_n^{I^c})_j$ with $\Swap^j_m(w) = \emptyset$. If $n$ is odd, there is at most one $t \in [n-1]$ for which $w$ has a $(-1,1)^j_t$ pair, and if such a pair exists, then $i(t+1) = 1$ and $\plussetblock_{1,i(t)-1} \in \opblock^j(w).$
\end{lem}
\begin{proof}
In Lemma ~\ref{pair_setup}, it was shown that there is at most one $t \in [n-1]$ for which $w$ has a $(-1,1)_t^j$ pair. All that remains to check is that $i(t+1) = 1,$ and $\plussetblock_{1,i(t)-1} \in \opblock^j(w).$
First, we will show that $i(t+1) = 1.$ By Lemma ~\ref{pair_setup}, we know $w_{1,j(1)} = 1.$ Note that $j(1)$ lies on an odd column right of $j$, which means $j(1) -1 \notin D(w)$ Hence, $w_{j(1)-1,i(j(1)-1)} = -1$, which implies that we have a $(-1,1)^j_{j(1)-1}$ pair. Since there is at most one such pair by Lemma ~\ref{pair_setup}, we must have $t = j(1)-1$. So, $i(t+1) = i(j(1))=1.$ 

The last claim to show is that if we do have a unique $(-1,1)^j_t$ pair, then $\plussetblock_{1,i(t)-1} \in \opblock^j(w).$ By Lemma ~\ref{pair_setup}, there must be some $p$ and $s$ with $\plussetblock_{1,p} \in \opblock^j(w)$ and $B_{1,s} \in \eblock^j(w),$ with $p < s.$ So, $p+1 \in B_{1,s}$ and thus $w_{p+1,j(p+1)} = -1.$ Since $p+1 \not\equiv n$ and all rows above $p+1$ have $1$s, $w$ has a $(-1,1)^j_{j(p+1)}$ pair. Since $w$ has only a single $(-1,1)^j_t$ pair, we must have $t = j(p+1),$ as claimed.
\end{proof}

\sssec{Property (5)}
\begin{lem}
\label{semi_adjacent_minus_ones}
Let $w\in (C_n^{I^c})_j$ with $\Swap^j_m(w) = \emptyset$. If there is $t \in [n-1]_0$ so that $w$ has a $(-1,-1)^j_t$ pair, then either $i(t+1) +1 = i(t)$ or $B_{i(t+1)+1,i(t)-1} \in \epblock^j(w).$
\end{lem}
\begin{proof}
First, we will show that there is some $B_{r,s} \in \block^j(w)$ with $i(t)$ and $i(t+1) \in B_{r,s}$. Second, we will show $B_{i(t+1)+1,i(t)-1} \in \epblock^j(w).$ The proofs of both of these facts involve an infinite descent.

Suppose there is no $r$ and $s$ so that $B_{r,s} \in \eblock^j(w)$ with $i(t)$ and $i(t+1) \in B_{r,s}$. Let $r$ and $s$ be so that $B_{r,s} \in \eblock^j(w)$ with $t \in B_{r,s}.$ By Lemma ~\ref{even_minus_ones}, $|\minus_{B_{r,s}} \cap s_e^j(w)| =|\minus_{B_{r,s}}\cap s_o^j(w)|.$ Using this, since $i(t+1) \notin B_{r,s},$ there must be some $t_2$ so that $w$ has a $(-1,-1)^j_{t_2}$ pair with $i(t_2+1) \in B_{r,s}$ but $i(t_2) \notin B_{r,s}.$ This implies there is some $B_{r_2,s_2} \in \block^j(w)$ with $s_2 > s$. However, replacing $t$ by $t_2$ and $B_{r,s}$ by $B_{r_2,s_2},$ we are in the same situation as we started in, but with $s_2 < s$. Repeating this process, we deduce that $\block^j(w)$ has unbounded size. Of course, this is a contradiction, as the matrix is finite.

At last, we only have to show $i(t+1) = i(t) -1$ or $B_{i(t+1)+1,i(t)-1} \in \epblock^j(w).$ Suppose $w$ has a $(-1,-1)^j_t$ pair with $i(t+1) \neq i(t)-1$ and $B_{i(t+1)+1,i(t)-1} \notin \epblock^j(w)$. By the preceding paragraph, we have some $B_{r,s} \in \eblock^j(w)$ with $i(t)$ and $i(t+1) \in B_{r,s}$. Let $k_1$ be the maximum row index so that $k_1$ is above $i(t)$ and $w_{k_1,j(k_1)} = -1$. 
Let $s_1 = j(k_1)$. Note that $s_1 \not\equiv t,$ since either $i(s_1) = i(t)-1$ or $B_{i(s_1)+1,i(t)-1} \in \epblock^j(w).$ Further, note that $s_1 \neq t+1$ by assumption that $i(t+1) \neq i(t)-1$ and $B_{i(t+1)+1,i(t)-1} \notin \epblock^j(w)$. Since $s_1 \equiv n,$ we have $s_1-1 \notin D(w),$ and $w$ has a $(-1,-1)^j_{s_1-1}$ pair. By construction, $i(s_1-1)<i(t)$ because $i(s_1-1) < i(s_1)$ and $i(s_1-1) \neq i(t).$ Hence, we are in the same situation we started in, but with $(t,t+1)$ replaced by $(s_1-1,s_1)$ and $i(s_1-1) <i(t).$ So, we have an infinite descent, a contradiction.
\end{proof}

\sssec{Property (6)}

\begin{lem}
\label{no_top_even_minus_block}
Let $w\in (C_n^{I^c})_j$ with $\Swap^j_m(w) = \emptyset$. If $B_{r,s} \in \block^j(w)$ with $r \neq 1,$ then there does not exist $\minussetblock_{p,q} \in \emblock^j(w)$ with $\minussetblock_{p,q} \subsetneq B_{r,s}$ and $p \equiv r.$
\end{lem}
\begin{proof}

Suppose $r >1$ and there is such a $\minussetblock_{p,q} \in \emblock^j(w),$ with $p \equiv r.$ Then, all $\plussetblock_{l,m} \in \pblock^j(w)$ with $\plussetblock_{l,m} \subset B_{r,s}$ satisfy $\plussetblock_{l,m} \in \epblock^j(w)$, by Lemma ~\ref{even_ones} and all $\minussetblock_{p,q} \in \mblock^j(w)$ with 
$p \neq r,q \neq s$ satisfy $\minussetblock_{p,q} \in \emblock^j(w)$ by Lemma~\ref{semi_adjacent_minus_ones}. These facts, together with $p \equiv r$ imply that all $\minussetblock_{p,q} \in \mblock^j(w)$ and $\minussetblock_{p,q}\subset B_{r,s}$ must satisfy $\minussetblock_{p,q} \in \emblock^j(w).$ Now, fix $p,q$ such that $\minussetblock_{p,q} \in \mblock^j(w)$ with $\minussetblock_{p,q} \subset B_{r,s}$ so that for any other $\minussetblock_{l,m} \in \mblock^j(w)$ with $\minussetblock_{l,m} \subset B_{r,s}$ we have $p > l$. There are now two cases depending on whether $p \equiv n.$

\vskip.1in
\noindent
{\sf Case 1: $p \equiv n$}

Since $\minussetblock_{p,q} \neq B_{r,s},$ either $p-1 \in B_{r,s},$ or $\plussetblock_{q+1,s} \in \epblock^j(w)$. However, we know that $p-1 \not \equiv n$ and $s \not \equiv n.$ If $w_{p-1,j(p-1)} = 1,$ since $w_{p,j(p)} = -1,$ we have a contradiction to Lemma ~\ref{adjacent_ones}. If $w_{s,j(s)} = 1,$ then $j(s+1)$ is left of $j,$ which again contradicts Lemma ~\ref{adjacent_ones}. Hence, this situation is impossible.

\vskip.1in
\noindent
{\sf Case 2: $p \not\equiv n$}

In this case, note that $q \equiv n.$ Hence $j(q)-1 \not \equiv n$ and $i(j(q)-1) > q.$ But this means that $i(j(q)-1) \notin B_{r,s},$ as $\minussetblock_{p,q}$ is the lowest $j$-block of $-1$s in $w$ contained in $B_{r,s}.$ Therefore, $w$ has a $(-1,-1)^j_{j(q)-1}$ pair, which is not contained in any element of $\eblock^j(w)$, contradicting Lemma ~\ref{semi_adjacent_minus_ones}.
\end{proof}

\begin{cor}
\label{minus_one_structure}
Let $w\in (C_n^{I^c})_j$ with $\Swap^j_m(w) = \emptyset$. Let $B_{r,s} \in \block^j(w)$ with $r \neq 1$ and $B_{p,q} \in \mblock^j(w)$. If $\minussetblock_{p,q} \subsetneq B_{r,s},$ then $\minussetblock_{p,q} \in \omblock^j(w)$ if and only if $p=r$ or $q = s.$
\end{cor}
\begin{proof}
By Lemma ~\ref{even_differences}, $B_{r,s} \in \eblock^j(w)$.
By Lemma ~\ref{no_top_even_minus_block}, if $p = r$ or $q = s,$ we must have that $\minussetblock_{p,q} \in \omblock^j(w)$.
Here we are using that if $q = s$ and $\minussetblock_{p,q} \in \emblock^j(w)$, since $B_{r,s} \in \eblock^j(w)$, we have $p \equiv r$. However, by Lemma ~\ref{no_odd_minus}, if $p \neq q$ and $r \neq s,$ then $\minussetblock_{p,q} \in \emblock^j(w)$.
\end{proof}

\sssec{Property (7)}

\begin{lem}
\label{no_top_even_ones}
Let $w\in (C_n^{I^c})_j$ with $\Swap^j_m(w) = \emptyset$. Let $r \neq 1,$ let $B_{r,s}\in \eblock^j(w),$ and let $\plussetblock_{r,q} \in \pblock^j(w)$ \textup{(}respectively $\plussetblock_{p,s} \in \pblock^j(w)$\textup{)}, then $q = s$ \textup{(}respectively $p = r$.\textup{)}
\end{lem}
\begin{proof}
Suppose there exists $\plussetblock_{r,q}\subset B_{r,s}$ with $\plussetblock_{r,q} \in \epblock^j(w)$, with $q \neq s.$ Then, there is some $t \in [n-1]$ so that $\minussetblock_{q+1,t} \in \emblock^j(w),$ with $q+1 \equiv r,$ contradicting Lemma ~\ref{no_top_even_minus_block}.

An analogous argument goes through in the case of $\plussetblock_{p,s} \in \epblock^j(w):$ We then have some $t \in [n-1]$ with $t \equiv r$ for which $\minussetblock_{t,s-1} \in \mblock^j(w),$ contradicting Lemma ~\ref{no_top_even_minus_block}.
\end{proof}
\sssec{Equivalence of the Seven Properties and an empty set of General-Minus $j$-Swaps}
\label{sssec:equivalence_gen_minus}

\begin{thm}
\label{equivalence_uncancelled}
An element $w \in (C_n^{I^c})_j$ is $j$-initially uncancelled if and only if $\Swap^j_m(w) = \emptyset.$
\end{thm}
\begin{proof}
First, let us check that if $\Swap^j_m(w) = \emptyset,$ then it is $j$-initially uncancelled.
Property (1) holds by Lemma ~\ref{n-1_row}. Property (2) holds by Lemma ~\ref{even_differences}. Property (3) holds by Lemma ~\ref{even_ones}. Property (4) is satisfied by Lemmas ~\ref{no_even_pairs} and ~\ref{one_odd_pair} . Property (5) follows from Lemma ~\ref{semi_adjacent_minus_ones}. By Corollary ~\ref{minus_one_structure}, Property (6) holds. Finally, by Lemma ~\ref{no_top_even_ones}, Property (7) holds.

The converse also follows in a straightforward fashion. Lemma ~\ref{equivalence_even_differences} shows that if $\oblock^j(w) = \emptyset$, then $\Swap^j_l(w) = \emptyset$. So, we only need to show $\Swap^j_\dplus(w) = \emptyset$, $\Swap^j_\dminus(w) = \emptyset,$ and $\Swap^j_\sminus(w) = \emptyset$. If $\Swap^j_\dplus(w) \cup \Swap^j_\dminus(w) \neq \emptyset$, then there must be some $B_{r,s} \in \block^j(w)$ with either $\plussetblock_{p,q} \in \opblock^j(w)$ or $\minussetblock_{p,q} \in \omblock^j(w)$ with $p \neq r, q \neq s.$ However, by Property (3), the former cannot happen. By Property (6), the latter cannot happen.

Finally, we check $ \Swap^j_\sminus(w) = \emptyset.$ Suppose $(b,t) \in  \Swap^j_\sminus(w)$ with $t \in B_{r,s}$. If $k \in B_{r,s}$ with $w_{k,j(k)} = -1,$ then by Properties (6) and (7), we must have both $w_{r,j(r)} =-1$ and $w_{s,j(s)} = -1.$ Therefore, we cannot have $(b,t) \in \Swap^j_\sminus$ with $i(t) \in B_{r,s},$ as this would imply either $w_{r,j(r)} = 1$ or $w_{s,j(s)} = 1.$
\end{proof}

\begin{cor}
\label{initially_uncancelled_support}
Recall $U_n^{I^c}$ is the set of $j$-initially uncancelled elements. It follows that
\begin{align*}
	\sum_{w \in (C_n^{I^c})_j} (-1)^{l(w)}X^{L(w)} = \sum_{w \in (U_n^{I^c})_j} (-1)^{l(w)}X^{L(w)}.
\end{align*}
\end{cor}
\begin{proof}
To prove this, we only need to show there is an involution $i:(C_n^{I^c})_j \setminus (U_n^{I^c})_j \rightarrow (C_n^{I^c})_j \setminus (U_n^{I^c})_j,$ which is sign reversing and $L$ preserving. However, we have shown in Theorem ~\ref{equivalence_uncancelled} that $w\in (C_n^{I^c})_j \setminus (U_n^{I^c})_j$ if and only if $\Swap_m^j(w) \neq \emptyset.$ Therefore, using Lemma \ref{restrict_swaps}, the involution $\involution^j_m$ defined in Lemma ~\ref{involution} restricts to a sign reversing, $L$ preserving involution on $(C_n^{I^c})_j \setminus (U_n^{I^c})_j$.
\end{proof}

\subsection{Finally Uncancelled Elements}

Analogously to initially uncancelled elements, we next define finally uncancelled elements.
Again, these are defined by $7$ opaque conditions, but turn out to precisely characterize those
elements of $C_n^{I^c}$ with no general-plus $k$-swaps, as shown in Theorem ~\ref{finally_uncancelled_support}

\begin{defn}
\label{finally_uncancelled}
Let $I = \{i_1,\ldots, i_l\}$. Fix $k \in I \cup \{n\}$ such that $k\equiv i$ for all $i \in I \cup \{n\}$ with $i>k.$ An element $w \in C_n^{I^c},$ is {\it $k$-finally uncancelled} if $w$ satisfies the following properties.
\begin{enumerate}
	\item Either $j(n)$ is left of $k$ or $w_{n,j(n)} = 1.$
	\item $\oblock^k(w) = \emptyset.$
	\item Suppose $t$ is a row index so that $w_{t,j(t)}=-1,$ $j(t)$ is right of $k,$ and there is no $s \in [n]$ with $s \geq t$ for which $\minussetblock_{1,s} \in \omblock^k(w)$. Then, there exists $\minussetblock_{r,s} \in \emblock^k(w)$ with $r \leq t \leq s.$
	\item If $n$ is even, there can be at most one $t \in [n-1]$ for which $w$ has a $(-1,1)^k_t$ pair, and if $w$ has such a pair, then $i(t) = 1$ and $\minussetblock_{1,i(t+1)-1} \in \mblock^k(w)$. If $n$ is odd, there is no $t$ for which $w$ has a $(-1,1)^k_t$ pair.
	\item If $w$ has a $(1,1)^k_t$ pair, then either $i(t)+1 = i(t+1)$ or $\minussetblock_{i(t)+1,i(t+1)-1} \in \mblock^k(w).$
	\item Suppose $\plussetblock_{p,q} \in \pblock^k(w)$ and $B_{r,s} \in \block^k(w)$ with $r>1$ and $\plussetblock_{p,q} \subsetneq B_{r,s}.$ Then $\plussetblock_{p,q} \in \omblock^j(w)$ is an odd block of $1$s if and only if $p=r$ or $q = s.$
	\item If $r >1$ and $s$ are row indices with $B_{r,s} \in \block^k(w)$, and $q$ is a row index with $\minussetblock_{r,q} \in \mblock^k(w)$ (respectively $\minussetblock_{p,s} \in \mblock^k(w)$), then $q = s$ (respectively $p = r$.)
\end{enumerate}
Notate $(V_n^{I^c})^k = \{w \in C_n^{I^c} \mid w \text{ is }k \text{-finally uncancelled }\}.$
\end{defn}

\begin{rem}
\label{initial_final_differences}
It is important to note that if $w \in (U_n^{I^c})_j,$ then $w(j) = n.$ However, for $w \in (V_n^{I^c})^k,$ we make no such restriction. The reason we were able to further assume that $w(j) = n$ when proving Theorem ~\ref{equivalence_uncancelled} is that the maximal element of $\Swap_m^j(C_n^{I^c})$ under the swapping ordering preserves column $j$. However, it will not always be the case that the maximal element of $\Swap_p^k(C_n^{I^c})$ preserves column $k$. That is, we have no analog of Lemma ~\ref{restrict_swaps}. Hence, in Theorem ~\ref{equivalence_uncancelled} we obtained that $(U^{I^c}_n)_j$ is precisely the set of $w \in (C_n^{I^c})_j$ with $\Swap^j_m(w) =\emptyset$, and $j(n)=1.$ In Theorem ~\ref{finally_uncancelled_support} we will obtain a ``weaker'' statement that $(V^{I^c}_n)^k$ is the set of $w \in (C_n^{I^c})$ satisfying $\Swap^k_p(w) = \emptyset$.
\end{rem}

\begin{thm}
\label{finally_uncancelled_support}
For $I$ and $k$ as in Definition \ref{finally_uncancelled}, 
\begin{align*}
	\sum_{w \in C_n^{I^c}} (-1)^{l(w)}X^{L(w)} = \sum_{w \in (V_n^{I^c})^k} (-1)^{l(w)}X^{L(w)}.
\end{align*}
\end{thm}
\begin{proof}
The proof of this result is essentially the same as that of Corollary ~\ref{initially_uncancelled_support}, but with the role of $-1$s and $1$s reversed. One key difference is that here we are using the set $\Swap_p^k(C_n^{I^c}),$ whereas Corollary ~\ref{initially_uncancelled_support} used $\Swap_m^j(C_n^{I^c}).$ See Remark ~\ref{initial_final_differences} for a further discussion of the differences. The analogs of Lemmas ~\ref{even_differences}, ~\ref{n-1_row}, ~\ref{even_ones}, ~\ref{no_even_pairs}, ~\ref{one_odd_pair}, ~\ref{semi_adjacent_minus_ones}, Corollary ~\ref{minus_one_structure}, and Lemma ~\ref{no_top_even_ones} all easily go through similarly with $1$s and $-1$s interchanged. Then, the analog of Theorem ~\ref{equivalence_uncancelled} tells us 
$\Swap^k_p(w)=\emptyset$, if and only if $w$ is $k$-finally uncancelled. Finally, we can use the involution $\involution^k_p$ to cancel off all elements of $C_n^{I^c} \setminus (V_n^{I^c})^k.$
\end{proof}

\ssec{The Lowering Map}
\label{subsection:lowering-map}

\begin{note}
\label{fixed_set_J}
Fix $I = \{i_1,\ldots,i_{r-1}, j,i_{r+1}\ldots, i_l\} \subset [n-1]_0,n \in I$ with $j = i_r$ and $\delta_k \equiv 0$ for $k \geq r$, as was stipulated in Notation ~\ref{fixed_j}. Then, denote
\begin{align*}
	J &= I^{(r)} =\{i_1,\ldots, i_{r-1},j-1,i_{r+1}-1,\ldots, i_l-1\} \subset [n-2]_0.
\end{align*}
\end{note}

{\bf For the rest of this section}, retain this notation for $I$ and $J.$ In order to complete the proof of Theorem ~\ref{induction_result}, it remains to show there is a bijection between $(U_n^{I^c})_j$ and $(V_{n-1}^{J^c})^{j-1}$ that decreases $L$ by $n-j$ and preserves the sign of $w.$ We define this map in Algorithm ~\ref{algorithm}, and then prove it satisfies the desired properties.

\begin{alg}
\label{algorithm}
The {\it lowering map} $H: (U_n^{I^c})_j\rightarrow B_{n-1}$ takes a $j$-initially uncancelled element $w \in (U_n^{I^c})_j$ to an element $H(w),$ where $H(w)$ is computed by the following algorithm:

\vskip.1in
\noindent
{\sf Step 1}

Let $w_1 = ws_j \cdots s_{n-2}s_{n-1}.$ Since $(w_1)_{n,n} = 1$, consider $w_1$ as in $B_{n-1},$ which is naturally embedded in $B_n$ as the set of all elements fixing $n.$

\vskip.1in
\noindent
{\sf Step 2}

For $1 \leq r < t \leq n$, define the function 
\begin{align*}
	p_{r,t} \colon B_{n-1} & \longrightarrow B_{n-1} \\
	x & \longmapsto s_rs_{r+1}\cdots s_tx.
\end{align*}
Let $\epblock^{j-1}(w_1) = \{\plussetblock_{r_1,s_1},\ldots, \plussetblock_{r_p,s_p}\}$ and $w_2 = p_{r_1,t_1} \circ \cdots \circ p_{r_p,s_p} (w_1).$
\vskip.1in
\noindent
{\sf Step 3}

If $\minussetblock_{1,1} \in \omblock^{j-1}(w_2)$, let $w_3 = s_0w_2.$ Otherwise, let $w_3 = w_2$.

\vskip.1in
\noindent
{\sf Step 4}

For $1 \leq r < t \leq n$, define the function 
\begin{align*}
	n_{r,t} \colon B_{n-1} & \longrightarrow B_{n-1} \\
	x & \longmapsto s_{t-1}\cdots s_{r-1}x.
\end{align*}
Let $\mblock^{j-1}(w_3) = \{B_{r_1,s_1},\ldots, B_{r_p,s_p}\}$ and $H(w)= n_{r_p,t_p} \circ \cdots \circ n_{r_1,s_1} (w_3).$
\end{alg}

\begin{rem}
Note that in Steps 2 and 4, $w_2$ and $H(w)$ are independent of the order in which we apply $p_{r,t},n_{r,t}$, since nonadjacent transpositions commute. 
\end{rem}
\begin{rem}
An alternative characterization of Algorithm ~\ref{algorithm} is as follows. First cross out column $j$ and row $n,$ to get a matrix in $B_{n-1}.$ Then, start at the bottom of the matrix, and move upward. Whenever there is a $1$ in row $t$ with $j(t)$ right of $j,$ switch rows $t$ and $t+1.$ Once we hit the top of the matrix, we proceed downward, and whenever row $t$ has a $-1$ with $j(t)$ right of $j,$ interchange rows $t$ and $t-1,$ unless $t = 1,$ in which case we change the sign of row $1$.

We tend to use this point of view of interchanging two adjacent rows in the ensuing proofs, although of course, it is equivalent to applying transpositions to the matrix.
\end{rem}

\begin{eg}
\label{alg_example}
Let us now give an example of how the algorithm would proceed with the matrix
\[w = \begin{tikzpicture}[baseline=-0.5ex]
\matrix(B)[matrix of math nodes,column sep={20pt,between origins},row
    sep={20pt,between origins},left delimiter={(},right delimiter={)}] at (6,0)
    {
	 &  &  &  &  &  & 1 &  & \\
	 &  &  &  &  & -1 &  &  & \\
	 1&  &  &  &  &  &  &  & \\
	 & 1 &  &  &  &  &  &  & \\
	 &  &  &  & -1 &  &  &  & \\
	 &  &  &  &  &  &  &  1& \\
	 &  &  &  &  &  &  &  & 1\\
	 &  &  & -1 &  &  &  &  & \\
	 &  & 1 &  &  &  &  &  & \\
	 };
  \end{tikzpicture}.
  \]
Here, since $w(3) = 9,$ we set $j = 3.$ In Step 1, we move $j(9)=3$ to the far right, resulting in
\[
	\begin{tikzpicture}[baseline=-0.5ex]
\matrix(B)[matrix of math nodes,column sep={20pt,between origins},row
    sep={20pt,between origins},left delimiter={(},right delimiter={)}] at (6,0)
    {
    &  &  &  &  & |[name=1l]|1 & |[name=1r]|\bullet &  & \\
	 &  &  &  & |[name=2l]|-1 & |[name=2r]|\bullet &  &  & \\
	 1&  &  &  &  &  &  &  & \\
	 & 1 &  &  &  &  &  &  & \\
	 &  &  & |[name=5l]|-1 & |[name=5r]|\bullet &  &  &  & \\
	 &  &  &  &  &  & |[name=6l]|1 & |[name=6r]|\bullet & \\
	 &  &  &  &  &  &  &  |[name=7l]|1& |[name=7r]|\bullet\\
	 &  & |[name=8l]|-1 & |[name=8r]|\bullet &  &  &  &  & \\
	 &  & |[name=9l]|\bullet&  &  &  &  &  & |[name=9r]|1\\
  };
  \draw[->]
  (1r)to(1l);
  \draw[->]
  (2r)to(2l);
  \draw[->]
  (5r)to(5l);
  \draw[->]
  (6r)to(6l);
  \draw[->]
  (7r)to(7l);
  \draw[->]
  (8r)to(8l);
  \draw[->]
  (9l)to(9r);
  
  \end{tikzpicture}.
  \]
Therefore, we obtain
\[w_1=\begin{tikzpicture}[baseline=-0.5ex]
\matrix(B)[matrix of math nodes,column sep={20pt,between origins},row
    sep={20pt,between origins},left delimiter={(},right delimiter={)}] at (6,0)
    {
	 &  &  &  &  & 1 &  &   \\
	 &  &  &  & -1 &  &  &  \\
	 1&  &  &  &  &  &  &   \\
	 & 1 &  &  &  &  &  &   \\
	 &  &  & -1 &  &  &  &  \\
	 &  &  &  &  &  & 1 &   \\
	 &  &  &  &  &  &  &  1\\
	 &  & -1 &  &  &  &  & \\
	 };
  \end{tikzpicture}.
  \]
Next, we perform Step 2, which moves the $(1,1)^2_7$ pair down by 1, and pushes the $-1$ in column 3 up by 2.
\[ 
	\begin{tikzpicture}[baseline=-0.5ex]
\matrix(B)[matrix of math nodes,column sep={20pt,between origins},row
    sep={20pt,between origins},left delimiter={(},right delimiter={)}] at (6,0)
    {
     &  &  &  &  & 1 &  &   \\
	 &  &  &  & -1 &  &  &  \\
	 1&  &  &  &  &  &  &   \\
	 & 1 &  &  &  &  &  &   \\
	 &  &  & -1 &  &  &  &  \\
	 &  & |[name=3u]|-1 &  &  &  & |[name=7u]|\bullet &   \\
	 &  &  &  &  &  & |[name=7d]|1 & |[name=8u]|\bullet\\
	 &  &  |[name=3d]|\bullet&  &  &  &  & |[name=8d]|1\\
  };
  \draw[->]
  (3d)to(3u);
  \draw[->]
  (7u)to(7d);
  \draw[->]
  (8u)to(8d);
  
  \end{tikzpicture}.
  \]
However, we have not completed Step 2, because there is another 
$2$-block of $1$s in 
$w_1$, namely $\plussetblock_{1,1} \in \pblock^2(w),$ and so we transpose the top two rows. The resulting matrix is
\[w_2 = 
	\begin{tikzpicture}[baseline=-0.5ex]
\matrix(B)[matrix of math nodes,column sep={20pt,between origins},row
    sep={20pt,between origins},left delimiter={(},right delimiter={)}] at (6,0)	
    {
    &  &  &  & |[name=5u]|-1 & |[name=6u]|\bullet&  &   \\
	 &  &  &  & |[name=5d]|\bullet &|[name=6d]| 1 &  &  \\
	 1&  &  &  &  &  &  &   \\
	 & 1 &  &  &  &  &  &   \\
	 &  &  & -1 &  &  &  &  \\
	 &  & -1 &  &  &  &  &   \\
	 &  &  &  &  &  & 1 & \\
	 &  &  &  &  &  &  & 1\\
  };
  \draw[->]
  (5d)to(5u);
  \draw[->]
  (6u)to(6d);  
  \end{tikzpicture}.
  \] 

Next, we apply Step 3. Since $\minussetblock_{1,1} \in \omblock^2(w)$, we change the sign of the top row. This results in the matrix
\[w_3 = 
	\begin{tikzpicture}[baseline=-0.5ex]
\matrix(B)[matrix of math nodes,column sep={20pt,between origins},row
    sep={20pt,between origins},left delimiter={(},right delimiter={)}] at (6,0)	
    {
    &  &  &  & |[name=5u]|1 & |[name=6u]|&  &   \\
	 &  &  &  & |[name=5d]| &|[name=6d]| 1 &  &  \\
	 1&  &  &  &  &  &  &   \\
	 & 1 &  &  &  &  &  &   \\
	 &  &  & -1 &  &  &  &  \\
	 &  & -1 &  &  &  &  &   \\
	 &  &  &  &  &  & 1 & \\
	 &  &  &  &  &  &  & 1\\
  };
   \draw[->] 
   (5u) circle (2mm);
  
  \end{tikzpicture}.
  \]
Finally, we perform Step 4, moving the $-1$s in rows 5 and 6 up by one, and moving the $1$ in the $4$th row down two. We obtain
\[H(w) = 
	\begin{tikzpicture}[baseline=-0.5ex]
\matrix(B)[matrix of math nodes,column sep={20pt,between origins},row
    sep={20pt,between origins},left delimiter={(},right delimiter={)}] at (6,0)
    {
   &  &  &  & 1 &  &  &   \\
	 &  &  &  &  & 1 &  &  \\
	 1&  &  &  &  &  &  &   \\
	 &  |[name=2u]|\bullet &  & |[name=4u]| -1 &  &  &  &   \\
	 &  &  |[name=3u]|-1 &  |[name=4d]|\bullet &  &  &  &  \\
	 &  |[name=2d]|1 &  |[name=3d]|\bullet &  &  &  &  &   \\
	 &  &  &  &  &  & 1 & \\
	 &  &  &  &  &  &  & 1\\
  };
  \draw[->]
  (3d)to(3u);
  \draw[->]
  (4d)to(4u);
   \draw[->] 
   (2u) to (2d);
  
  \end{tikzpicture}.
  \]

This is our final result. Lo and behold, $H(w)$ is a chessboard element with a very similar descent set to our original matrix, and furthermore, it is actually a $(j-1)$-finally uncancelled matrix. We can also see by direct computation that $\sgn(w) = \sgn(H(w)) = -1,$ and $L(w) = 15,L(H(w)) = 9.$ Since $n = 9$ and $j = 3,$ we have that $L(w)-L(H(w)) = 15 - 9 = 6 = n-j.$
This is no coincidence, and this lowering map defines a bijection between the $j$-initially uncancelled elements of $(C_n^{I^c})_j$ because the $(j-1)$-finally uncancelled elements of $C_{n-1}^{I^c}$, in a way that preserves $\sgn(w)$ and changes $L(w)$ by $n-j.$
\end{eg}

\begin{lem}
\label{even_minus_one_blocks}
Let $w \in (U_n^{I^c})_j.$ Then, the matrix $w_3,$ after Step 3 of Algorithm ~\ref{algorithm}, has $\omblock^{j-1}(w_3) = \emptyset.$ Furthermore, each $t$
for which $w_3$ has a $(-1,-1)^j_t$ pair satisfies $i(t+1) = i(t)-1.$ That is, the two $-1$s in adjacent columns also lie in adjacent rows.
\end{lem}
\begin{proof}
If some $B_{p,q} \in \block^{j-1}(w) \setminus (\emblock^{j-1}(w) \cup \epblock^{j-1}(w))$ and $p>1,$ we know by Properties (6) and (7) of $j$-initially uncancelled matrices that the pattern of $1$s and $-1$s contained in $B_{p,q},$ is as follows: At the top, there is a odd $(j-1)$-block of $-1$s in $w$. Beneath it, there are alternating even $(j-1)$-blocks of $1$s in $w$ and even $(j-1)$-blocks of $-1$s in $w$, until finally there is an odd $(j-1)$-block of $-1$s in $w$. Hence, moving all the even $(j-1)$-blocks of $1$s in $w_1$ down by a row in Step 2 forces the top row of each $(j-1)$-block of $-1$s in $w_3$ to move to a higher $(j-1)$-block of $-1$s in $w_3$. This, together with possibly changing the sign of the top row, as is done in Step $3,$ ensures that all of the resulting $(j-1)$-blocks of $-1$s in $w_3$ will lie in $\emblock^{j-1}(w_3)$. For an example, see the movement between $w_1$ and $w_3$ in Example ~\ref{alg_example}.

The fact that $i_{w_3}(t+1) = i_{w_3}(t)-1,$ follows directly from Property (5). That is, since we move the $1$s between each $(-1,-1)^j_{t+1}$ pair in $w$ down, the $(-1,-1)^j_{t+1}$ pair in $w$ becomes a $(-1,-1)^{j-1}_t$ pair in $w_3,$ which is moved together in Step 2 so that $i_{w_3}(t+1) = i_{w_3}(t)-1.$
\end{proof}

\begin{lem}
\label{descent_preservation}
The image $H(U_n^{I^c})_j$ is contained in $C_{n-1}^{J^c}.$
\end{lem}
\begin{proof}

We must show that for any $w \in (U_n^{I^c})_j,$ we have 
$D(H(w)) \subset J,$
and that the final result is a chessboard element. Fix some $w \in (U_n^{I^c})_j,$ and let $w_i$ denote the result after Step $i$ of Algorithm ~\ref{algorithm} is applied to $w$. First, it is clear that $D(w_1) \subset J,$ because $D(w) \subset I,$ and in Step 1 we have $t \in D(w_1)$ and $t > j$ if and only if $t+1 \in D(w).$ In Steps 2--4, the descent set also remains unchanged: Any column left of $j$ maintains the same relative position to any other column left of $j,$ which shows that $D(w) \cap [j-1]_0 = D(H(w)) \cap [j-1]_0$. For rows $r$ with $j(r)$ right of $j,$ we are only moving rows which contain $1$s up, and rows which contain $-1$s down, and we are always preserving the relative order of the rows with $-1$s and preserving the relative order of the rows with $1$s. So, the final descent set is contained in $J$ and the image of $H$ is contained in $B_{n-1}^{J^c}.$

Next, we check that the image only contains chessboard elements. After Step 1, all columns $t$ left of $j$ are on squares with $t \equiv i_{w_1}(t) $. For this argument, we will call squares $(i(t),t)$ with $t \equiv i(t)$ chessboard squares. After Step 1, all the $\pm 1$s which are right of column $j$ are not on chessboard squares.

In Step 2, we analyze the effect of $p_{r,s}$ on $w_1$ in two cases, depending on whether $\opblock^{j-1}(w_1) = \emptyset$.

First, consider the case that $\opblock^{j-1}(w_1) = \emptyset$. For each $B_{r,s} \in \epblock^{j-1}(w_1)$, we have that $p_{r,s}$ moves each $1$ in a column right of $j$ down by a single row, so that it will then be on a chessboard square. Since $B_{r,s} \in \epblock^{j-1}(w_1),$ for all $k \notin B_{r,s},(k,j(k)),$ we have that $k$ is chessboard before applying $p_{r,s}$ if and only if it is chessboard after applying $p_{r,s}$.

Second, consider the case that $\opblock^{j-1}(w_1) \neq \emptyset,$ so that $\{B_{1,s}\} = \opblock^{j-1}(w_1)$ by Property (3) of Definition ~\ref{initially_uncancelled_defn}. By applying $p_{1,s},$ the $-1$ in row $s+1$ is moved up an odd number of rows, and hence for all $m \in [s+1]$, we have $i_{w_2}(m) \not \equiv i_{w_1}(m).$

In Step 3, we see for all $k \in [n]$ that $i_{w_2}(k) = i_{w_3}(k)$.

Finally, in Step 4, if there is some $k$ with $\minussetblock_{1,k} \in \emblock^{j-1}(w_3)$, the $-1$ in the top row is changed to a 1 and moved down to the even row $k+1$. The other rows $m$ with $2 \leq m \leq s$ satisfy $i_{H(w)}(m) = i_{w_3}(m) +1$. Otherwise, if $\minussetblock_{p,q} \in \mblock^{j-1}(w_3),$ then the parity of all rows $m$ with $p \leq m \leq q$ satisfy $i_{H(w)}(m) = i_{w_3}(m) +1$. Finally, thanks to Lemma ~\ref{even_minus_one_blocks}, $\omblock^{j-1}(w_3) = \emptyset$, and so for $m$ not in any $\minussetblock_{p,q} \subset \oblock^{j-1}(w_3),$ we have $i_{H(w)}(m) \equiv i_{w_3}(m).$

Hence, from Steps $2$, $3$, and $4$, we have changed the parity of the rows $r$ with $j(r)$ right of $j$ and preserved the parity of the rows $r$ with $i(r) \leq j$, resulting in a chessboard element.
\end{proof}

\begin{lem}
\label{correct_image}
The image $H((U_{n}^{I^c})_{j})$ is contained in $(V_{n-1}^{J^c})^{j-1}$.
\end{lem}
\begin{proof}Let $w \in (U_{n}^{I^c})_{j}$ and let $w_i$ be the result of Algorithm ~\ref{algorithm} after Step $i$. 
We check that $H(w)$ satisfies the seven properties of $(j-1)$-finally uncancelled matrices, as defined in
Definition~\ref{finally_uncancelled}.
 
\vskip.1in
\noindent
{\sf Property (1)}

First, we show $H(w)$ satisfies Property (1). If $w_3$ does not satisfy Property (1), then in Step 4, the $-1$ in row $n-1$ is moved up by a row, and so either $j(n-1)\leq j-1$ or $H(w)_{n-1,j(n-1)}=1.$

\vskip.1in
\noindent
{\sf Property (2)}

Next, $H(w)$ satisfies Property (2) because the parities of the rows of nonzero entries in columns left of $j$ are preserved by Algorithm ~\ref{algorithm}, as shown in Lemma ~\ref{descent_preservation}. Further, the relative positions of the rows $r$ with $j(r)$ left of $j$ are also preserved under Algorithm ~\ref{algorithm}. 
It follows that the parities of the number of rows in blocks in $H(w)$ agree with those in $w$. In particular, if there are no odd $j$-blocks in $w$, there will also be no odd $(j-1)$-blocks in $H(w)$.

\vskip.1in
\noindent
{\sf Property (3)}

Property (3) holds by Lemma ~\ref{even_minus_one_blocks}, since Step 4 of the algorithm does not separate any $j$-blocks of $-1$s in $w$.

\vskip.1in
\noindent
{\sf Property (4)}

In order to show Property (4) of $(j-1)$-finally uncancelled matrices holds after the algorithm, we break up the argument into five cases depending on what the topmost $j$-block $B_{r,s} \in \block^j(w)$ is. That is, we can check it in the following five cases: (a) $r > 1,$ (b) there is some $q$ with $\minussetblock_{1,q} \in \omblock^j(w),$ (c) there is some $q$ with $\minussetblock_{1,q} \in \opblock^j(w),$ (d) there is some $q$ with $\minussetblock_{1,q} \in \emblock^j(w),$ and finally (e) $\minussetblock_{1,q} \in \epblock^j(w).$ Below, we will only write out the details for case (d). The other cases are similar.

Suppose there is some $q$ with $\minussetblock_{1,q} \in \emblock^j(w)$. Since row $1$ does not lie in an odd $j$-block of $-1$s in $w$, there are no $(-1,1)^j_t$ pairs. So, $w$ has a $(-1,-1)^j_{j_w(1)-1}$ pair, meaning $1 \equiv n$. Then, we see $\minussetblock_{1,q-1} \in \oblock^{j-1}(H(w))$ and $H(w)$ has a $(-1,1)^{j-1}_{j_{H(w)}(1)}$ pair. Additionally, $i_{H(w)}(j_{H(w)}(1)) = 1,$ (as is always true,) and $i_{H(w)}(j_{H(w)}(1)+1) = q+1,$ satisfying Property (4) of $(j-1)$-finally uncancelled matrices.

\vskip.1in
\noindent
{\sf Property (5)}

Property (5) holds because if we start with $w \in B_n^{I^c}$ which has a $(1,1)^j_t$ pair, then by Lemma ~\ref{adjacent_ones} we initially have $i_w(t+1) = i_w(t)+1.$ Therefore, these 1s will remain adjacent through Step 3, that is, $i_{w_3}(t) = i_{w_3}(t-1)+1.$ Therefore, in Step 4, they will either remain adjacent, or be separated by a single $(j-1)$-even block of $-1$s in $H(w)$.

\vskip.1in
\noindent
{\sf Properties (6) and (7)}

It is easy to check these conditions are satisfied when we are dealing with $B_{r,s} \in \block^j(w)$ with $r = 1.$ For the remainder of this verification, we assume $r > 1.$

Observe next that Properties (6) and (7) of $(j-1)$-finally uncancelled matrices are equivalent to the following: For $B_{r,s} \in \block^{j-1}(H(w))$, which contains both $1$s and $-1$s,
\begin{align}
\label{block_decomposition}
B_{r,s} = \plussetblock_{r,s_1} \cup  \minussetblock_{s_1+1,s_2}\cup \plussetblock_{s_2+1,s_3} \cup \cdots \cup \minussetblock_{s_{k-2}+1,s_{k-1}}\cup \plussetblock_{s_{k-1}+1,s},
\end{align}
where $\plussetblock_{r,s_1},\plussetblock_{s_{k-1}+1,s} \in \opblock^{j-1}(H(w)),$ but $\plussetblock_{s_2+1,s_3},\plussetblock_{s_4+1,s_5}, \ldots, \plussetblock_{s_{k-3}+1, s_{k-2}} \in \epblock^{j-1}(H(w))$ and 
$\minussetblock_{s_1+1,s_2},\minussetblock_{s_3+1,s_4},\ldots,\minussetblock_{s_{k-2}+1,s_{k-1}} \in \emblock^{j-1}(w).$

By Properties (2), (3), and (5), to complete the proof, we only have to show that if $B_{r,s} \in \block^{j-1}(H(w))$ is given a decomposition as in \eqref{block_decomposition} then $\plussetblock_{r,s_1} \notin \epblock^{j-1}(H(w))$. Suppose there is some $w \in (U_n^{I^c})_j$ so that this is the case. Working backwards, we see $\plussetblock_{r,s_1-1} \in \opblock^{j-1}(w_3).$ There are now two cases: First, if $r +2 = s_1$, we see there are some $t,k,l \in [n]$ with $B_{r,t} \in \block^j(w)$ with $B_{r,k} \in \emblock^j(w)$ and $B_{k,l} \in \epblock^j(w).$ But this means $w$ does not satisfy Property (6) of $j$-initially uncancelled matrices.
Second, if $s_1 - r > 2,$ we obtain that there are some $t,k,l$ so that $B_{r,t} \in \block^j(w), B_{r,s_1-2} \in \epblock^j(w)$ and $B_{s_1-2,s_2-2} \in \emblock^j(w),$ implying that $w$ does not satisfy Property (7) of $j$-initially uncancelled matrices.
\end{proof}

\ssec{The Lowering Map's Effect on $L$}

In this subsection, we check that $H$ preserves $\sgn(w)$ and changes $L(w)$ by $n-j.$

\begin{lem}
Let $w \in (U_n^{I^c})_j$. Then, $\sgn(w) = \sgn(H(w))$. 
\end{lem}
\begin{proof}
In all four steps, we multiply $w$ by a certain number of Coxeter generators. It suffices to check that we multiply $w$ by an even number of Coxeter generators to obtain $H(w)$. In Step 1, there are $n-j\equiv 0$ Coxeter generators by which $w$ is multiplied to obtain $w_1.$ In Steps $2$--$4$, $w$ is also multiplied by $n-j$ Coxeter generators. Hence, $w$ is multiplied by a total of $2(n-j)$ Coxeter generators to obtain $H(w),$ and hence $\sgn(w) = \sgn(H(w))$.
\end{proof}

\begin{lem}
\label{step_1_calculation}
Recall the notation from
Definition~\ref{us_defn}.
Let $w \in (C_n^{I^c})_j,$ and $y = s_{n-1}s_{n-2}\cdots s_{j}w.$ Note that $y(n)=n.$ Then, $L(y) = L(w)-|u_o^j(w)s_e^j(w)| - |u_e^j(w)s_o^j(w)| + |u_e^j(w)s_e^j(w)| + |u_o^j(w)s_o^j(w)| - \frac{n-j}{2}.$
\end{lem}

\begin{proof}

We calculate 
\begin{align*}
	a(y) -a(w) &= -|s_o^j(w)|+|s_e^j(w)|,\\ 
	b(y)-b(w) &= -\frac{n-j}{2},\\
	c(y)-c(w)&=-|u_o^j(w)s_e^j(w)| + |u_o^j(w)s_o^j(w)|.
\end{align*}

Using the abc Lemma ~\ref{abc_lemma}, $L(w) = a(w) + b(w) + 2c(w),$ so 
\begin{align*}
	L(y) - L(w)&=a(y) +b(y) + 2c(y) - a(w) - b(w) - 2c(w) \\
	& = a(y) - a(w) + b(y) - b(w) + 2(c(y) - c(w))
	\\
	&=-|s_o^j(w)|+|s_e^j(w)|-\frac{n-j}{2} -2|u_o^j(w)s_e^j(w)| + 2 |u_o^j(w)s_o^j(w)| 
	\\
	&= -|u_o^j(w)s_e^j(w)| - |u_e^j(w)s_o^j(w)| + |u_e^j(w)s_e^j(w)| + |u_o^j(w)s_o^j(w)| - \frac{n-j}{2}. \qedhere
\end{align*}
\end{proof}

\begin{lem}
\label{local_change}
Let $\epsilon \in \left\{ \pm 1 \right\}$. Then,
\begin{equation}
\label{local_1}
	L\begin{pmatrix}
	0 & 1\\
	\epsilon & 0
\end{pmatrix}
-L\begin{pmatrix}
	\epsilon & 0\\
	0 & 1
\end{pmatrix}=1
\end{equation}
\begin{equation}
\label{local_2}
	L\begin{pmatrix}
	0 & -1\\
	\epsilon & 0
\end{pmatrix}
-L\begin{pmatrix}
	\epsilon & 0\\
	0 & -1
\end{pmatrix}=-1
\end{equation}

\begin{equation}
\label{local_3}
	L\begin{pmatrix}
	0 & 1\\
	-1 & 0
\end{pmatrix}
-L\begin{pmatrix}
	1 & 0\\
	0 & 1
\end{pmatrix}=2
\end{equation}

\begin{equation}
\label{local_4}
	L\begin{pmatrix}
	0 & -1\\
	-1 & 0
\end{pmatrix}
-L\begin{pmatrix}
	-1 & 0\\
	0 & 1
\end{pmatrix}=1.
\end{equation}
\end{lem}
\begin{proof}
The proof is given by calculating $L$ in the four situations.
\end{proof}

\begin{lem}
\label{adjacent_swaps}
Let $w \in (U_n^{I^c})_j.$ Let rows $t,t+1,t+2$ be three consecutive positive integers with $t \equiv n,$ and suppose there is some $q \geq t+2$ with $\plussetblock_{t+1,q} \in \pblock^j(w)$, \textup{(}respectively there is some $p \leq t$ with $\plussetblock_{p,t+1} \in \pblock^j(w)$\textup{)}. If we let $v = s_{t+1}s_tw$ (respectively $v = s_t s_{t+1}w$), then, $L(w)-1 = L(v)$.
\end{lem}
\begin{proof}
Let us prove this for the case where there is some $q \geq t+2$ with $\plussetblock_{t+1,q} \in \pblock^j(w)$ using Lemma ~\ref{adjacent_ones}. The other case in which there is $p \leq t$ with $\plussetblock_{p,t+1} \in \pblock^j(w)$ is analogous, by reversing the role of $1$s and $-1$s.

Since $t+2 \equiv n,$ it follows that $w$ has a $(1,1)^j_{j(t+1)}$ pair, and hence $j(t+1) = j(t)+1$. Define $w_t$ to be the $3 \times 3$ sub-matrix defined by rows $t,t+1,t+2$ and columns $j(t),j(t+1),j(t+2).$ Let $v_t$ be the analogous $3 \times 3$ sub-matrix of $v$.

First, since rows $t,t+1,t+2$ are adjacent, $L(w) - L(v)=L(w_t) - L(v_t).$ One can prove this using methods analogous to those in proof of the Swapping Lemma, by noting that the values of $b_j(w_t)=b_j(v_t)$ and $c_j(w_t)=c_j(v_t)$, as long as $j \notin \{j(t),j(t+1),j(t+2)\}.$

Now, there are two cases, either $j(t+2)<j(t)$ or $j(t+2)>j(t).$ For simplicity, let us assume $j(t+2)>j(t);$ the other case is analogous. By the preceding paragraph, we have $L(v)-L(w) = L(v_t)-L(w_t).$ We may note that $a(v_t) = a(w_t)$ and also that
\begin{align*}
	b_{2,3}(w_t) &= b_{2,3}(v_t) &\text{ and } &&c_{2,3}(w_t) = c_{2,3}(v_t).
\end{align*}
The only difference in $L(v_t)$ and $L(w_t)$ comes from columns $1$ and $2,$ which is precisely calculated by \eqref{local_1} and \eqref{local_2} from Lemma ~\ref{local_change}. Thus, $L(v_t) - L(w_t) = -1,$ which implies $L(v) - L(w) = -1.$
\end{proof}

\begin{lem}
\label{us_constraints}
For any $w\in (C_n^{I^c})_j$ with $\Swap^j_m(w) = \emptyset$, we have 
\begin{align*}
|u_o^j(w)s_e^j(w)| - |u_o^j(w)s_o^j(w)| &=0, \\ 
|u_e^j(w)s_e^j(w)| - |u_e^j(w)s_o^j(w)| &=
\begin{cases}
	1, &\text{ if there is some } t \text{ for which } w \text{ has a }(-1,1)^j_t \text{ pair},\\
	0 &\text{ otherwise.}
\end{cases}
\end{align*}

\end{lem}
\begin{proof}
Observe that by Lemma ~\ref{semi_adjacent_minus_ones}, $|u_o^j(w)s_o^j(w)| - |u_o^j(w)s_e^j(w)|$ is, up to sign, the number of $t$ for which $w$ has a $(-1,1)^j_t$ pair with $t \in u_o^j(w).$ Similarly, $|u_e^j(w)s_e^j(w)| - |u_e^j(w)s_o^j(w)|$ is, up to sign, the number of $(-1,1)^j_t$ pairs with $t \in u_e^j(w)$. By Lemmas ~\ref{one_odd_pair} and ~\ref{no_even_pairs} the only way we can have a $(-1,1)^j_t$ pair is if $n$ is odd, $t$ is even, and there are no rows $s$ with $s <i(t)$.
Therefore, if $w$ does have a $(-1,1)^j_t$ pair, then we must have $t \in u_e^j(w)$ and $t \in s_e^j(w).$ Hence, $|u_o^j(w)s_o^j(w)| - |u_o^j(w)s_e^j(w)|=0$ and $|u_e^j(w)s_e^j(w)| - |u_e^j(w)s_o^j(w)| =1$ precisely when there is some $t$ for which $w$ has a $(-1,1)^j_t$ pair.
\end{proof}

\begin{lem}
\label{L_change}
Let $w \in (U_n^{I^c})_j.$ Then, $L(H(w)) = L(w)-n+j.$
\end{lem}
\begin{proof}
By Lemma ~\ref{step_1_calculation}, if $w_1$ is the result of $w$ after Step 1 of Algorithm \ref{algorithm}, then 
\begin{align*}
L(w_1) = L(w)-|u_o^j(w)s_e^j(w)| - |u_e^j(w)s_o^j(w)| + |u_e^j(w)s_e^j(w)| + |u_o^j(w)s_o^j(w)| - \frac{n-j}{2}.
\end{align*}
So, we only need to show that in Steps 2--4, we change $L$ by 
\begin{align*}
|u_o^j(w)s_e^j(w)| + |u_e^j(w)s_o^j(w)| - |u_e^j(w)s_e^j(w)| - |u_o^j(w)s_o^j(w)| - \frac{n-j}{2}.
\end{align*}

We now proceed by casework, depending on whether or not there is some $t$ so that $w$ has a $(-1,1)^j_t$ pair.

\vskip.1in
\noindent
{\sf Case 1: There is no $(-1,1)^j_t$ pair}

Since there is no $t \in [n-1]$ so that $w$ has a $(-1,1)^j_t$ pair, by Lemma ~\ref{us_constraints}, \begin{align}\label{no_pair_even_odd}
|u_o^j(w)s_e^j(w)| + |u_e^j(w)s_o^j(w)| - |u_e^j(w)s_e^j(w)| - |u_o^j(w)s_o^j(w)|=0.
\end{align}
\newline
\noindent
{\sf Case 1a: $w_{1,j(1)}= 1.$}
In this case, we will make exactly $\frac{n-j}{2}$ moves of the type from Lemma ~\ref{adjacent_swaps}. By Lemma ~\ref{adjacent_swaps}, each such move decreases $L$ by 1. By \eqref{no_pair_even_odd}, $L$ changes by 
\begin{align*}
-\frac{n-j}{2}= |u_o^j(w)s_e^j(w)| + |u_e^j(w)s_o^j(w)| - |u_e^j(w)s_e^j(w)| - |u_o^j(w)s_o^j(w)|- \frac{n-j}{2}.
\end{align*}
\newline
\noindent
{\sf Case 1b: $w_{1,j(1)} = -1.$}

In this case, we will only have $\frac{n-j}{2}-1$ moves as in Lemma ~\ref{adjacent_swaps}. However, during Step 2, we will have one move as in \eqref{local_4}. Since this decreases $L$ by one, Steps 2--4 change $L$ by
\begin{align*}
\left(-\frac{n-j}{2}+1 \right	)-1=|u_o^j(w)s_e^j(w)| + |u_e^j(w)s_o^j(w)| - |u_e^j(w)s_e^j(w)| - |u_o^j(w)s_o^j(w)|- \frac{n-j}{2}.
\end{align*}

\vskip.1in
\noindent
{\sf Case 2: $w$ has a $(-1,1)^j_t$ pair}

In this case, we will have $\frac{n-j}{2}-1$ moves as in Lemma ~\ref{adjacent_swaps}. However, in Step 2, we have a single move as in \eqref{local_3}. Note that this move decreases $L$ by 2. Therefore, $L$ changes by $-\frac{n-j}{2}+1-2 = -\frac{n-j}{2}-1$, and indeed, 
\begin{align*}
|u_o^j(w)s_e^j(w)| + |u_e^j(w)s_o^j(w)| - |u_e^j(w)s_e^j(w)| - |u_o^j(w)s_o^j(w)|=-1.
\end{align*}
Therefore, in Steps 2 through 4, using Lemma ~\ref{us_constraints}, $L$ changes by 
\begin{align*}-\frac{n-j}{2}-1 = |u_o^j(w)s_e^j(w)| + |u_e^j(w)s_o^j(w)| - |u_e^j(w)s_e^j(w)| - |u_o^j(w)s_o^j(w)| - \frac{n-j}{2}.
\end{align*}
\end{proof}

\subsection{The Inverse to the Lowering Map}

We have shown so far that there is a map $H(U_n^{I^c})_j \rightarrow (V_{n-1}^{J^c})^{j-1}$ that alters $L(w)$ by $n-j$ and preserves $\sgn(w).$ In order to complete the proof of Theorem ~\ref{induction_result}, we only need to show this is a bijection. In what follows, we produce the inverse map.

\begin{alg}
Define the {\it raising map} $R: (V_{n-1}^{J^c})^{j-1}\rightarrow B_n,$ which takes a $(j-1)$-finally uncancelled element $w \in (V_n^{J^c})^{j-1}$, to an element $R(w),$ where $R(w)$ is computed by the following algorithm.

\vskip.1in
\noindent
{\sf Step 1}

For $1 \leq r < t \leq n$, define the functions 
\begin{align*}
	rn_{r,t}\colon B_{n-1} & \longrightarrow B_{n-1}\\
	x & \longmapsto s_r\cdots s_tx.
\end{align*}
Let $\emblock^{j-1}(w) = \{\minussetblock_{r_1,s_1},\ldots, \minussetblock_{r_p,s_p}\}.$ Let $w_1= rn_{r_1,t_1} \circ \cdots \circ rn_{r_p,s_p} (w).$

\vskip.1in
\noindent
{\sf Step 2}

If $\plussetblock_{1,1} \in \opblock^{j-1}(w_1)$, let $w_2 = s_{0}w_1.$

\vskip.1in
\noindent
{\sf Step 3}

For $1 \leq r < t \leq n$, define the functions 
\begin{align*}
	rp_{r,t} \colon B_{n-1} & \longrightarrow B_{n-1}\\
	x & \longmapsto s_{t-1} \cdots s_{r-1}x.
\end{align*}
Let $\pblock^{j-1}(w_2) = \{\plussetblock_{r_1,s_1},\ldots, \plussetblock_{r_p,s_p}\}$ and $w_3= rp_{r_p,t_p} \circ \cdots \circ rp_{r_1,s_1} (w_2).$

\vskip.1in
\noindent
{\sf Step 4}

Let $w_4$ be the image of $w_3$ under the natural embedding $B_{n-1} \rightarrow B_n$. That is,
$$w_4(i) =\begin{cases}
	w_3(i), &\text{ if }i\neq \pm n,\\
	i &\text{ otherwise.}
\end{cases}$$

\vskip.1in
\noindent
{\sf Step 5}

Define $L(w) = w_4 s_{n-1}\cdots s_{j+1}s_{j}.$
\end{alg}

\begin{lem}
\label{correct_preimage}
The image $R((V_{n-1}^{J^c})^{j-1})$ is contained in $(U_{n}^{I^c})_{j}.$
\end{lem}
\begin{proof}
The proof is analogous to that of Lemma ~\ref{correct_image}.
\end{proof}
\begin{lem}
\label{bijection}
The maps $H:(U_n^{I^c})_j\rightarrow(V_{n-1}^{J^c})^{j-1}, R:(V_{n-1}^{J^c})^{j-1}\rightarrow (U_n^{I^c})_j$ are mutual inverses.
\end{lem}
\begin{proof}
By Lemma ~\ref{correct_image} and Lemma ~\ref{correct_preimage}, $H$ is indeed a map from $(U_n^{I^c})_j$ to $(V_{n-1}^{J^c})^{j-1}$ and $R$ is indeed a map from $(V_{n-1}^{J^c})^{j-1}$ to $(U_n^{I^c})_j.$ Furthermore, the algorithms were precisely constructed so that $H \circ R = \id$ and $R\circ H= \id,$ since in $H$ we are multiplying by transpositions $s_i,$ and in $R$ we are multiplying by the same transpositions in the opposite order. Here we are crucially using that $s_i^2 = 1,$ as $B_n$ is a Coxeter group.
\end{proof}

\ssec{Completion of Part 3}

\begin{cor}
\label{part_3_induction_result}
Part 3 of Theorem ~\ref{induction_result} holds. That is, let $I \subset [n-1]_0,$
and $j = i_k \in I \cup \{n\}.$ Then $S_{n,I,j}(X)= X^{n-j} \cdot f_{n-1,J}(X),$ where $J$ is as defined in Notation ~\ref{fixed_set_J}.
\end{cor}
\begin{proof}
We proceed by induction. Note that Theorem ~\ref{induction_result} trivially holds when $n = 1$, so we will assume it holds for $n-1$ and prove it holds for $n$.

In Lemma ~\ref{bijection}, we showed that $H:(U_n^{I^c})_j\rightarrow(V_{n-1}^{J^c})^{j-1}$ is a bijection.
By Lemma ~\ref{L_change}, we see $L(H(w)) = L(w) -n +j$ and $\sgn(H(w)) = \sgn(w).$
This tells us that
\begin{equation}
\begin{split}
\label{initial_to_final}
	\sum_{w \in (U_n^{I^c})_j} (-1)^{l(w)}X^{L(w)} = \sum_{w \in (U_n^{I^c})_{j}} (-1)^{l(H(w))}X^{L(H(w))+n-j} \\
	= X^{n-j}\cdot \sum_{w \in (V_{n-1}^{J^c})^{j-1}} (-1)^{l(w)}X^{L(w)}.
\end{split}
\end{equation}
Therefore,
 \begin{align*}
	 S_{n,I,j}(X) &= \sum_{w \in (C_n^{I^c})_j} (-1)^{l(w)}X^{L(w)}  & (\text{by Definition} ~\ref{snik_defn})
	\\
	&= \sum_{w \in (U_n^{I^c})_j} (-1)^{l(w)}X^{L(w)} & (\text{by Corollary} ~\ref{initially_uncancelled_support})
	\\
	&= X^{n-j}\cdot \sum_{w \in (V_{n-1}^{J^c})^{j-1}} (-1)^{l(w)}X^{L(w)} & (\text{by } ~\eqref{initial_to_final})
	\\
	&= X^{n-j} \cdot\sum_{w \in C_{n-1}^{J^c}} (-1)^{l(w)}X^{L(w)} &(\text{by Theorem} ~\ref{finally_uncancelled_support})
	\\
	&= X^{n-j} \cdot f_{n-1,J}(X) & (\text{by the inductive hypothesis}).
\end{align*}
\end{proof}

\section{A Further Conjecture}
\label{further_conjecture}

In the above sections, we found a nice way to factor $\sum_{w \in B_n^{I^c}} (-1)^{l(w)}X^{L(w)}.$ It is natural to ask if there is any nice way to factor the expression $\sum_{w \in B_n^{I^c}} t^{l(w)}X^{L(w)},$ where we replace $-1$ by a general variable $t.$ It seems that in general, there is not a nice factorization, but in the case that $0 \in I,$ the expression does have a nice factor.

\begin{prop}
\label{one_direction_general_t_conj}
If $0 \in I$, then the two variable polynomial $Xt+1$ divides $\sum_{w \in B_n^{I^c}} t^{l(w)}X^{L(w)}$.
\end{prop}
\begin{proof}
The following proof is due to an anonymous referee.
Assume $0 \in I$.
Let $P_n^{I^c} \subset B_n^{I^c}$ denote those elements $w \in B_n^{I^c}$ with $w_{1, j(1)} = 1$ and let
$M_n^{I^c} \subset B_n^{I^c}$ denote those elements $w \in B_n^{I^c}$ with $w_{1, j(1)} = -1$
(where $P$ stands for plus and $M$ stands for minus).
Note that $B_n^{I^c} = P_n^{I^c} \coprod M_n^{I_c}$, since every $w \in B_n^{I_c}$ either has a $1$
or $-1$ in its first row.
Because of our assumption that $0 \in I$, left multiplication by the Coxeter element $s_0$ defines an involution
\begin{align*}
	\iota \colon B_n^{I^c} & \rightarrow B_n^{I^c} \\
	w & \mapsto s_0w.
\end{align*}

We now record several pleasant properties of the involution $\iota$.
First, by definition, $\iota(P_n^{I^c}) = M_n^{I^c}$ and $\iota(M_n^{I^c}) = P_n^{I^c}$.
Next, for $w \in P_n^{I^c}$, we have $L(s_0 w) = L(w) + 1$.
One can see this directly from the definition of $L$ or by noting that $a(s_0 w) = a(w) + 1, b(s_0w) = b(w),$ and $c(s_0w) = c(w)$.
Similarly, for $w \in P_n^{I^c}$, $l(s_0w) = l(w) + 1$.
This follows, for example, from the explicit description of $l$ given in 
\cite[Proposition 8.1.1 and (8.2)]{bjorner}:
\begin{align*}
	l(w) = \left | \left\{ \left( i,j \right) \in \left[ n \right]\times \left[ n \right] \mid i < j, w(i) > w(j) \right\} \right | 
	+ \left | \left\{ \left( i,j \right) \in \left[ n \right]\times \left[ n \right] | i \leq j, w(-i) > w(j) \right\} \right |.
\end{align*}

Putting the above together, we find
\begin{align*}
	\sum_{w \in B_n^{I^c}} t^{l(w)}X^{L(w)} &=\sum_{w \in P_n^{I^c}} t^{l(w)}X^{L(w)} + \sum_{w \in M_n^{I^c}} t^{l(w)}X^{L(w)} \\
	&= \sum_{w \in P_n^{I^c}} \left( t^{l(w)}X^{L(w)} + t^{l(s_0 w)}X^{L(s_0 w)} \right) \\
	&= \sum_{w \in P_n^{I^c}} \left( t^{l(w)}X^{L(w)} + t^{l(w) + 1}X^{L(w)+1} \right) \\
	&= \sum_{w \in P_n^{I^c}} \left( t^{l(w)}X^{L(w)} + tX \cdot t^{l(w)}X^{L(w)} \right) \\
	&= (1+tX) \cdot \sum_{w \in P_n^{I^c}} t^{l(w)}X^{L(w)}.\qedhere
\end{align*} \end{proof}

It is natural to ask whether the converse of Proposition ~\ref{one_direction_general_t_conj} holds,
and we conjecture that it does.

\begin{conj}
\label{general_t_conj}
The two variable polynomial $Xt+1$ divides $\sum_{w \in B_n^{I^c}} t^{l(w)}X^{L(w)}$ if and only if $0 \in I.$
\end{conj}

\begin{rem}
Conjecture ~\ref{general_t_conj} is true for all subsets $I \subset [n-1]_0$ with $n \leq 6$, as was verified by a computer. However, we do not see a way to generalize the techniques used in this paper to prove Conjecture ~\ref{general_t_conj}. Note also that when we plug in $t=-1,$ we obtain the formula in Theorem ~\ref{result}, but $-X+1$ divides $f_{n,I}(X)$ so long as $I \neq \emptyset.$ As such, Theorem ~\ref{result} does not appear to be very helpful in proving Conjecture ~\ref{general_t_conj}.
\end{rem}
\begin{rem}
	One potential approach to proving Conjecture ~\ref{general_t_conj} would be to use the map given by switching the sign of the first row,
	described in the proof
	of Proposition ~\ref{one_direction_general_t_conj}. This does not define an involution (or even an automorphism) of $B_n^{I^c}$ when $0 \notin I$, but
	it does define an involution when restricted to those elements with $w_{1,1} \neq 1$. So, one could reduce the support of the sum
	to those elements with $w_{1,1} = 1$. However, it is not clear how to proceed from here.
\end{rem}

\section{Acknowledgements}

I would like to thank Vic Reiner for bringing the conjecture to my attention, offering many helpful edits, and for listening to my ideas. Of course, I thank Alexander Stasinski, Christopher Voll, and Rob Little for the many hours they spent scrutinously reading over this proof. I deeply thank several anonymous referees for helping me clarify and streamline numerous confusing issues and notation, and for catching many errors in the paper. 
I also thank a referee for pointing out the proof of Proposition ~\ref{one_direction_general_t_conj}.
I thank David Hemminger and Zijian Yao for carefully critiquing the most technical parts of the argument. Further thanks goes to my father, Peter Landesman, for his useful edits and his support. Finally, I would like to thank Levent Alpoge for an especially motivating comment regarding starting work before the REU. This research was conducted at the University of Minnesota, Twin Cities REU, supported by RTG grant NSF/DMS-1148634 and the Univ. of Minnesota School of Mathematics.

\bibliographystyle{abbrv}
\bibliography{biblio}

\end{document}